\documentclass[11pt,onecolumn]{article}
\usepackage{color,graphicx,xcolor,cite}
\usepackage{amsfonts,amssymb,amsmath,amsthm}
\usepackage{empheq}
\usepackage[ruled]{algorithm2e}
\usepackage{url}
\usepackage{multirow,array,booktabs}
\usepackage{threeparttable}
\usepackage{cases} 

\usepackage{latexsym}
\usepackage{color}
\definecolor{myblue}{rgb}{0,0,0.5}
\definecolor{mygreen}{rgb}{0,0.5,0}
\definecolor{myred}{rgb}{0.5,0,0}

\usepackage{amsmath,amsthm,amssymb}
\usepackage{graphicx}
\usepackage{hyperref}
\newcommand{\nn}{\nonumber}

\def \[{\begin{equation}}
\def \]{\end{equation}}

\newtheorem{theorem}{Theorem}[section]

\newtheorem{lemma}{Lemma}[section]

\newtheorem{remark}{Remark}[section]

\def \Q{{\mbox{${\cal Q}$}}}
\def \pd{{\mbox{\footnotesize$_{\!{P\!D}}$}}}
\def \dup{{\mbox{\footnotesize$_{\!{D\!P}}$}}}

\def \M{{\mbox{${\cal M}$}}}

\def \G{{\mbox{${\cal G}$}}}

\usepackage{ulem}

\begin{document}

\begin{center}

{\LARGE  \bf On construction of splitting contraction algorithms in a prediction-correction framework for separable convex optimization}\\

\bigskip
\medskip

  {\bf Bingsheng He}\footnote{\parbox[t]{15cm}{
 Department of Mathematics,  Nanjing University, Nanjing, China.
  This author was supported by the NSFC Grant 11871029. Email: hebma@nju.edu.cn}}
  \quad and \quad
 {\bf Xiaoming Yuan}\footnote{\parbox[t]{16cm}{
 Department of Mathematics, The University of Hong Kong, Hong Kong. Email:  xmyuan@hku.hk
  }}

\medskip

 \today

\end{center}

\medskip

{\small

\parbox{0.95\hsize}{

{\bf Abstract.} In the past decade, we had developed a series of splitting contraction algorithms for separable convex optimization problems, at the root of the alternating direction method of multipliers. Convergence of these algorithms was studied under specific model-tailored conditions, while these conditions can be conceptually abstracted as two generic conditions when these algorithms are all unified as a prediction-correction framework. In this paper, in turn, we showcase a constructive way for specifying the generic convergence-guaranteeing conditions, via which new splitting contraction algorithms can be generated automatically. It becomes possible to design more application-tailored splitting contraction algorithms by specifying the prediction-correction framework, while proving their convergence is a routine.

\medskip


 \medskip

  }}

\section{Introduction}

\setcounter{equation}{0}

To understand the convergence rate of the well-known alternating direction method of multipliers (ADMM) proposed in \cite{GM}, we had initiated an analytic technique in \cite{HeYuanSIAMNM} and further used it to conduct convergence analysis for a series of ADMM-based splitting contraction algorithms for separable convex programming problems in our previous works such as \cite{HeLWY,HeTaoYuanSIAMO,HeTaoYuanIMA,HeYuanNMD,HeXuYuan}. As summarized in \cite{HeGXJS16} (see also \cite{He20Years,HeYuanCOA18}), the key of this analytic technique is to represent the ADMM or its variants in \cite{HeLWY,HeTaoYuanSIAMO,HeTaoYuanIMA,HeYuanNMD,HeXuYuan} by a prediction-correction framework in the context of the variational inequality reformulation of the convex programming model under discussion, and then to analyze the contraction property of the iterative scheme to derive the convergence. In particular, various specific conditions for ensuring convergence of the algorithms in \cite{HeLWY,HeTaoYuanSIAMO,HeTaoYuanIMA,HeYuanNMD,HeXuYuan} can be conceptually abstracted as two generic conditions. Since the ADMM is an application of the Douglas-Rachford splitting method in \cite{DRSM,LM} and convergence of the mentioned algorithms can be proved by the contraction property, we call them ``splitting contraction algorithms" uniformly. Our previous result in \cite{CHYY} showing the divergence of the direct extension of the ADMM to multiple block separable convex programming problems was also based on this analytical technique.

In this paper, we revisit the unified framework in \cite{HeGXJS16} (see also \cite{He20Years,HeYuanCOA18,HeXuYuan}), including the prediction-correction representation of various splitting contraction algorithms and the generic convergence-guaranteeing conditions, and showcase how to constructively specify the generic convergence-guaranteeing conditions. Through the procedure, new splitting contraction algorithms can be automatically generated via specifying the prediction-correction framework. The constructive way provides flexibility to designing various application-tailored splitting contraction algorithms for given specific separable convex optimization models, while proving their convergence is a routine.

\section{Unified framework}

\setcounter{equation}{0}

We consider the canonical convex minimization problem with linear constraints
   \[  \label{Problem-LC}
    \min \{ \theta(u) \; | \;  {\cal A}u=b, \; u\in {\cal U}\}
  \]
where $\theta:\Re^{n}\to {\Re}$ is a closed, proper, and convex but not necessarily smooth function, ${\cal U}\subseteq\Re^n$ is a closed convex set, ${\cal {A}}\in\Re^{m\times n}$, and $b\in\Re^m$. For further analysis, we recall the prediction-correction representation of splitting contraction algorithms in the existing literatures \cite{HeYuanCOA18,HeLWY,HeTaoYuanSIAMO,HeTaoYuanIMA,HeXuYuan} in the context of variational inequalities, and the generic convergence-guaranteeing conditions.

\subsection{Variational inequality reformulation}

The Lagrangian function of \eqref{Problem-LC}  is
\[ \label{Lagran-A}
       L(u,\lambda) = \theta(u) - \lambda^T({\cal A}u-b)
       \]
with $\lambda \in \Re^m$ the Lagrange multiplier. A pair of  $(u^*,\lambda^*)\in {\cal U}\times \Re^m$ is called
  a saddle point of $L(u,\lambda)$ if
$$    L_{\lambda\in\Re^m}(u^*,\lambda)\le  L(u^*,\lambda^*)\le L_{u\in {\cal U}} (u,\lambda^*). $$
As analyzed in, e.g. \cite{He20Years,HeYuanCOA18,HeLWY,HeYuanSIAMNM,HeTaoYuanSIAMO,HeTaoYuanIMA,HeYuanNMD}, $(u^*,\lambda^*)$  can be characterized as a solution point of the following variational inequality (VI):
\[  \label{VIu-O}
      w^*\in \Omega, \quad \theta(u) -\theta(u^*) + (w-w^*)^T F(w^*) \ge 0, \quad \forall \,  w\in
      \Omega.
     \]
where
\[ \label{Notation-uFO}
    w = \biggl(\begin{array}{c}
                     u\\
                  \lambda \end{array} \biggr),
  \quad
    F(w) =\biggl(\begin{array}{c}
     - {\cal A}^T\lambda \\
      {\cal A}u-b \end{array} \biggr) \quad
  \hbox{and} \quad      \Omega= {\cal U} \times \Re^m.
   \]
We denote by $\Omega^*$ the set of solution points of the VI (\ref{VIu-O})-(\ref{Notation-uFO}).

\subsection{Prediction-correction framework and convergence conditions}

In the context of the VI (\ref{VIu-O})-(\ref{Notation-uFO}), the splitting contraction algorithms proposed in \cite{HeYuanCOA18,HeLWY,HeTaoYuanSIAMO,HeTaoYuanIMA,HeXuYuan} can all be uniformly represented by a prediction-correction framework, and their specific convergence-guaranteeing conditions can be conceptually abstracted as two generic conditions. Below, note that $v=w$, or $v$ could be a subvector of $w$, and ${\cal V}^*$ is denoted as the set of all $v^*$ whose corresponding $w^*\in \Omega^*$. Dimensionality of the matrices $Q$ and $M$ is conformed with that of $v$.

\begin{center}\fbox{\begin{minipage}{15.5cm}

\begin{subequations} \label{Framework}
\noindent{\bf \normalsize Prediction-correction framework for the VI (\ref{VIu-O})-(\ref{Notation-uFO})}

[Prediction Step.] With given $v^k$, find a vector $\tilde{w}^k\in \Omega$ such that
 \[ \label{Frame-Q}
    \theta(u) - \theta(\tilde{u}^k) +  (w- \tilde{w}^k )^T F(\tilde{w}^k)  \ge
      (v - \tilde{v}^k )^T Q (v^k-\tilde{v}^k),  \;  \forall  w\in\Omega,\]
where the matrix $Q$ is not necessarily symmetric but $Q^T +Q$ is assumed to be positive definite.

\smallskip
[Correction Step.] Find a nonsingular matrix $M$ and update $v$ by
 \[ \label{Frame-M}
     v^{k+1} = v^k -
      M (v^k-\tilde{v}^k).
   \]
\end{subequations}

\end{minipage}
}
\end{center}

\begin{center}\fbox{\begin{minipage}{15.5cm}
\noindent{\bf \normalsize Convergence conditions}

For the matrices $Q$ and $M$ used in (\ref{Frame-Q}) and (\ref{Frame-M}), respectively, there exists a matrix $H\succ0$ such that
 \begin{subequations} \label{Frame-HG}
\[   \label{Frame-H} H M=Q,\]
and
\[   \label{Frame-G}
G:=  Q^T +Q - M^THM \succ 0.
       \]
\end{subequations}
\end{minipage}
}
\end{center}

\begin{remark}
Note that the matrices $Q$ and $M$ take specific forms for the splitting contraction algorithms in \cite{HeYuanCOA18,HeLWY,HeTaoYuanSIAMO,HeTaoYuanIMA,HeXuYuan} when different separable cases of the model (\ref{Problem-LC}) are studied. With the specifically given $Q$ and $M$, it is easy to verify whether or not the conditions (\ref{Frame-HG}) are satisfied. How to construct the matrices $H$ and $M$ with a given $Q$ to ensure the convergence conditions (\ref{Frame-HG}), however, has not yet been discussed in the literature. This reverse procedure leads to the specification of a concrete splitting contraction algorithm for the model (\ref{Problem-LC}) via the prediction-correction framework (\ref{Framework}). A trivial case is when $Q$ is known to be even symmetric and positive definite. Then, we can simple take $M=\alpha I$ with $\alpha \in(0,2)$ and then $H= \frac{1}{\alpha} Q$, which certainly ensures the conditions \eqref{Frame-HG}. This case is indeed an application of the classic proximal point algorithm in \cite{Mar70}, as studied in \cite{CaiGuHeYuan,GuHY,HeGXJS22,HeYuanSIAM-IMS}. We would emphasize that the corresponding matrix $Q$ is generally asymmetric when many algorithms including those in \cite{HeYuanCOA18,HeLWY,HeTaoYuanSIAMO,HeTaoYuanIMA,HeXuYuan} are represented by the prediction-correction framework (\ref{Framework}). Hence, we focus on the nontrivial case where $Q$ is asymmetric to discuss how to constructively specify the prediction-correction framework (\ref{Framework}) and the convergence conditions (\ref{Frame-HG}).
\end{remark}

\subsection{Convergence}

Convergence of the prediction-correction framework (\ref{Framework}) under the conditions (\ref{Frame-HG}) has been shown in our previous works such as \cite{HeGXJS16,He20Years,HeXuYuan}. The main contraction property ensuring the convergence is summarized in the following theorem; we include the proof for completeness.

\begin{theorem} \label{THM-HauptA} Let $\{v^k\}$ be the sequence generated by the prediction-correction framework (\ref{Framework}) under the conditions (\ref{Frame-HG}). Then, it holds that
\begin{equation} \label{Haupt-Convergence} \|v^{k+1}-v^*\|_H^2
       \le \|v^k -v^*\|_H^2 -  \|v^k-\tilde{v}^k\|_G^2, \quad \forall v^*\in {\cal V}^*.
\end{equation}
\end{theorem}
\noindent{\bf Proof}. Using $Q=HM$ (see \eqref{Frame-H}), the prediction step can be written as
$$\tilde{w}^k\in \Omega, \quad   \theta(u) - \theta(\tilde{u}^k) +  (w- \tilde{w}^k )^T F(\tilde{w}^k)  \ge
      (v - \tilde{v}^k )^T HM (v^k-\tilde{v}^k),  \quad   \forall  w\in\Omega.  $$
Then, it follows from \eqref{Frame-M} that
$$
  \tilde{w}^k\in \Omega, \quad    \theta(u) - \theta(\tilde{u}^k) +  (w- \tilde{w}^k )^T\! F(\tilde{w}^k)
      \ge  (v- \tilde{v}^k)^T\!H(v^k -
    v^{k+1}),  \quad \forall w\in \Omega. $$
Setting $w=w^*$ in the above inequality, we get
  $$
   (v^k -  v^{k+1})^T\! H(\tilde{v}^k-v^*)\ge  \theta(\tilde{u}^k)-\theta(u^*) +  (\tilde{w}^k-w^* )^T\! F(\tilde{w}^k),
      \quad  \forall w^*\in \Omega^*. $$
Note that $(\tilde{w}^k-w^* )^T\! F(\tilde{w}^k) =(\tilde{w}^k-w^* )^T\! F(w^*)$ and the optimality of $w^*$. Thus,
we have
\[   \label{Tvv*} (v^k -
    v^{k+1})^TH (\tilde{v}^k-v^*)      \ge 0, \quad \forall  v^* \in {\cal V}^*.  \]
Setting   $a=v^k, \,  b=v^{k+1},    c=\tilde{v}^k$ and  $d=v^*$  in the identity
$$
   2 (a-b)^TH(c-d) = \bigl\{\|a-d\|_H^2 -\|b-d\|_H^2\bigr\} -
    \bigl\{ \|a-c\|_H^2 -\|b-c\|_H^2 \bigr\},
$$
we know from \eqref{Tvv*} that
\[ \label{LEM-M1-1} \|v^k -v^*\|_H^2- \|v^{k+1}-v^*\|_H^2  \ge  \|v^k-\tilde{v}^k\|_H^2 - \|v^{k+1}-\tilde{v}^k\|_H^2. \]
For the right-hand side of the last inequality, we have
\begin{eqnarray}  \label{LEM-M1-3}
\lefteqn{ \|v^k-\tilde{v}^k\|_H^2 - \|v^{k+1}-\tilde{v}^k\|_H^2}
   \nonumber \\
   &=&  \|v^k-\tilde{v}^k\|_H^2 - \|(v^k-\tilde{v}^k) -(v^k-v^{k+1})\|_H^2
   \nonumber \\
  & \stackrel{\eqref{Frame-M}}{=}  &  \|v^k-\tilde{v}^k\|_H^2 - \|(v^k-\tilde{v}^k)
  -   M (v^k-\tilde{v}^k)\|_H^2  \nonumber \\
  & = & 2(v^k-\tilde{v}^k)^TH M (v^k-\tilde{v}^k) - (v^k-\tilde{v}^k)^T
  M ^THM (v^k-\tilde{v}^k) \nonumber \\
  & = &  (v^k-\tilde{v}^k)^T ( Q^T +Q -  M^THM)(v^k-\tilde{v}^k) \nonumber \\
    &  \stackrel{\eqref{Frame-G}}{=} &  \|v^k-\tilde{v}^k\|_G^2.
\end{eqnarray}
Substituting  \eqref{LEM-M1-3}
  in  \eqref{LEM-M1-1}, the assertion \eqref{Haupt-Convergence} is proved. \qquad {\normalsize$\Box$}

\subsection{Naming of matrices}

Because of clear reasons, we call $Q$ and $M$ the prediction matrix and correction matrix, respectively. Moreover, the inequality \eqref{Haupt-Convergence} compares the proximity of two consecutive iterates to $v^*$ in the $H$-norm, and measures their difference (or progress) by the term $\|v^k-\tilde{v}^k\|_G^2$. Hence, we call $H$ and $G$ the norm matrix and profit matrix, respectively.

\subsection{Example}

To see how the prediction-correction framework \eqref{Framework} and the convergence conditions \eqref{Frame-HG} are related to our previous works
such as \cite{CHYY,HeLWY,HeTaoYuanSIAMO,HeTaoYuanIMA,HeXuYuan,HeYuanSIAMNM,HeYuanCOA18,HeYuanNMD}, we take the strictly contractive Peaceman-Rachford splitting method (PRSM) proposed in \cite{HeLWY} as an illustrative example.

\subsubsection{SC-PRSM}

We consider a two-block separable convex programming model
\[  \label{SCCO}
    \min \{\theta_1(x) + \theta_2(y) \; | \;  Ax + By =b, \, x\in {\cal X}, y\in{\cal Y} \},
   \]
where $\theta_i: {\Re}^{n_i}\to {\Re} \;(i=1,2)$ are closed, proper, and convex but not necessarily smooth functions; ${\cal X} \subseteq \Re^{n_1} $ and ${\cal Y} \subseteq \Re^{n_2}$ are closed convex sets; $A \in \Re^{m\times n_1}$ and $B \in \Re^{m\times n_2}$; $b\in \Re^m$; and $n=n_1+n_2$. This is a separable case of  \eqref{Problem-LC} with
   $$  u= \left(\begin{array}{c}
                     x\\  y \end{array} \right),  \quad  \theta(u) = \theta_1(x) + \theta_2(y),
                     \quad {\cal A} = (A, B) \quad \hbox{and}\quad  {\cal U}={\cal X}\times{\cal Y}. $$
Thus, the VI (\ref{VIu-O})-(\ref{Notation-uFO}) can be specified as the following:
  \begin{subequations} \label{VI-FORM-2}
\[  \label{VI-S2}
      w^*\in \Omega, \quad \theta(u) -\theta(u^*) + (w-w^*)^T F(w^*) \ge 0, \quad \forall \,  w\in
      \Omega.
     \]
where
 \[  w = \left(\begin{array}{c}
                     x\\  y\\  \lambda \end{array} \right),  \quad
                        F(w) = \left( \!\begin{array}{c}
            - A^T\lambda \\
             -B^T\lambda \\
             Ax + By -b \end{array} \! \right) \quad  \hbox{and} \quad
              \Omega={\cal X} \times {\cal Y}\times \Re^m.
             \]
\end{subequations}

Let the augmented Lagrangian function of \eqref{SCCO} be
\[ \label{A-Lagrange-2}
 {\cal L}^{\hbox{\scriptsize[2]}}_\beta(x,y,\lambda) =\theta_1(x) + \theta_2(y) - \lambda^T(Ax
 +By-b) +\frac{\beta}{2}     \|Ax+By-b \|^2,
 \]
with $\lambda\in \Re^m$ the Lagrange multiplier and $\beta >0$ the penalty parameter.
The well-known alternating direction method of multipliers (ADMM) proposed by Glowinski and Marrocco in \cite{GM} for solving \eqref{SCCO} is
\[  \label{ADMM-2}
\left\{ \begin{array}{l}
  x^{k+1}\in \arg\min \{{\cal L}^{\hbox{\scriptsize[2]}}_{\beta}(x,y^k,\lambda^k)\,|\, x\in {\cal X}\},\\[0.1cm]
  y^{k+1} \in \arg\min \{{\cal L}^{\hbox{\scriptsize[2]}}_{\beta}(x^{k+1},y,\lambda^{k}) \,|\, y\in {\cal Y}\},\\[0.1cm]
   \lambda^{k+1} =  \lambda^k - \beta (Ax^{k+1} + By^{k+1} -b).
   \end{array}\right.
   \]
As shown in \cite{Gabay}, the ADMM (\ref{ADMM-2}) is an application of the Douglas-Rachford splitting method in \cite{DRSM,LM} to the dual of (\ref{SCCO}). In \cite{HeLWY}, we proposed the following strictly contractive PRSM (SC-PRSM):

\begin{subequations} \label{S-ADMM-xlyl}
\begin{numcases}{\hskip-1.5cm \hbox{ (SC-PRSM)}\quad}
\label{S-ADMM-x} x^{k+1}\in \arg\min\{{\cal L}^{\hbox{\scriptsize[2]}}_{\beta}(x,y^k,\lambda^k)\,|\, x\in {\cal X}\},\\
\label{S-ADMM-l1} \lambda^{k+\frac{1}{2}}=\lambda^k- \mu\beta (Ax^{k+1}+By^k-b),\\
\label{S-ADMM-y}y^{k+1} \in \arg\min\{{\cal L}^{\hbox{\scriptsize[2]}}_{\beta}(x^{k+1},y, \lambda^{k+\frac{1}{2}}) \,|\, y\in {\cal Y}\},\\
\label{S-ADMM-l2}{\lambda}^{k+1} =   \lambda^{k+\frac{1}{2}}- \mu \beta (Ax^{k+1} + By^{k+1} -b),
\end{numcases}
\end{subequations}
where $\mu\in (0,1)$. It has been shown that if $\mu=1$ in (\ref{S-ADMM-xlyl}), then the resulting scheme (which can be regarded as a symmetric version of the ADMM (\ref{ADMM-2})) is an application of the PRSM in \cite{PRSM} to the dual of (\ref{SCCO}) while its convergence is not guaranteed as shown in \cite{CY}. For the SC-PRSM (\ref{S-ADMM-xlyl}) with $\mu \in (0,1)$, it was shown in \cite{HeLWY} that its sequence of iterations is strictly contractive with respect to the solution set and hence its convergence is guaranteed.

\subsubsection{Prediction-correction representation}

The follow lemma is useful and its proof can be found in, e.g. \cite{Beck}.

\begin{lemma} \label{CP-TF}
\begin{subequations} 
 Let ${\cal X}\subset \Re^n$ be a closed convex set, $\theta(x)$ and $f(x)$ be convex functions and $f(x)$ is
  differentiable  on an open set which contains ${\cal X}$. Assume that the solution set of the minimization problem
  $\min \{ \theta(x) + f(x) \, |\, x\in {\cal X}\}$ is nonempty.
     Then, {\color{black}
  \[  \label{CP-TF1} x^*  \in   \arg \min \{ \theta(x) + f(x) \, |\, x\in {\cal X}\}
     \]}
if and only if
\[\label{CP-TF2}
      x^*\in {\cal X}, \; \; \theta(x) -\theta(x^*) + (x-x^*)^T\nabla f(x^*) \ge 0, \; \; \forall\, x\in {\cal X}.
      \]\end{subequations}
\end{lemma}

Let us define $\tilde{w}^k$ as
\begin{equation} \label{Tnotation}
 \tilde{w}^k=
 \left(\begin{array}{c}
         \tilde{x}^k \\[0.1cm]
         \tilde{y}^k \\[0.1cm]
         \tilde{\lambda}^k
       \end{array}\right)
   = \left(
     \begin{array}{c}
        x^{k+1}\\[0.1cm]
        y^{k+1} \\[0.1cm]
        \lambda^k-\beta(Ax^{k+1}+By^{k}-b)
 \end{array}
  \right),
\end{equation}
where $(x^{k+1},y^{k+1})$ is generated by \eqref{S-ADMM-xlyl}. Now, we show that the SC-PRSM (\ref{S-ADMM-xlyl}) with $\tilde{w}^k$ defined in (\ref{Tnotation}) can be represented by the prediction-correction framework (\ref{Framework}).

\begin{theorem} \label{THM-HauptXLYL}  The SC-PRSM (\ref{S-ADMM-xlyl}) corresponds to the prediction-correction framework (\ref{Framework}) with
\[  \label{MatrixTQ}
      Q=  \left(\begin{array}{cc}
                \beta B^TB& -\mu B^T\\
         -B  & \frac{1}{\beta} I_m \end{array} \right) \quad \hbox{and}\quad M  =  \left(\begin{array}{cc}
              I    &    0 \\
             -\mu\beta B   &   2\mu I_m \end{array} \right).\]
\end{theorem}

\noindent {\bf Proof.} It follows from (\ref{Tnotation}) that $\tilde{\lambda}^k=\lambda^k-\beta(Ax^{k+1}+By^{k}-b)$. Also, the optimality condition of the $x$-subproblem \eqref{S-ADMM-x} is
\[  \label{ADMM-xlyl-x}
    \theta_1(x) - \theta_1(\tilde{x}^k)
       + (x - \tilde{x}^k)^T(-A^T\tilde{\lambda}^k) \ge 0, \;\; \forall\;  x\in {\cal X}.  \]
Notice that the objective function of the $y$-subproblem  \eqref{S-ADMM-y} is
$$
  {\cal L}^{[2]}_{\beta}(\tilde{x}^k,y,{\lambda}^{k+\frac{1}{2}})
    =  \theta_1(\tilde{x}^k) +\theta_2(y) -
    ({\lambda}^{k+\frac{1}{2}})^T(A\tilde{x}^k+By-b) +{\textstyle{\frac{\beta}{2}}}\|A\tilde{x}^k+By-b \|^2.
      $$
Ignoring some constant term in the $y$-subproblem, we obtain
$$  \tilde{y}^k=\arg\min\{\theta_2(y) -  ({\lambda}^{k+\frac{1}{2}})^T By + {\textstyle{\frac{\beta}{2}}}\|A\tilde{x}^k+By-b \|^2
        \,|\, y\in {\cal Y}\}.
         $$
Consequently, according to Lemma  \ref{CP-TF}, we have
$$
  \tilde{y}^k\in {\cal Y} , \quad  \theta_2(y) - \theta_2(\tilde{y}^k)
                 +  (y - \tilde{y}^k)^T  \bigl\{  -B^T{\lambda}^{k+\frac{1}{2}}  + \beta B^T(A\tilde{x}^k + B\tilde{y}^k-b)  \bigr\} \ge 0, \; \forall\, y\in {\cal Y}.
      $$
Using $\tilde{\lambda}^k=\lambda^k-\beta(Ax^{k+1}+By^{k}-b)$, we get
$$
     \lambda^{k+\frac{1}{2}}={\lambda}^k -\mu(\lambda^k-\tilde{\lambda}^k)=\tilde{\lambda}^k +(1- \mu)(\lambda^k-\tilde{\lambda}^k),
 $$
and
$$ \beta(A\tilde{x}^k + By^k-b)=  (\tilde{\lambda}^k -\lambda^k).   $$
Finally,  the optimality condition of the $y$-subproblem can be written as  $\tilde{y}^k\in {\cal Y}$, and
\[  \label{ADMM-xlyl-y}
   \theta_2(y) - \theta_2(\tilde{y}^k) + (y - \tilde{y}^k)^T \bigl\{ -B^T{\tilde{\lambda}}^k + \beta B^TB(\tilde{y}^k - y^k)
             - \mu B^T(\tilde{\lambda}^k -\lambda^k) \bigr\}\ge 0, \;\; \forall y\in {\cal Y}.
     \]
According to the definition of $\tilde{w}^k$ in \eqref{Tnotation}, and $\tilde{\lambda}^k=\lambda^k-\beta(Ax^{k+1}+By^{k}-b)$, we have
\[  \label{ADMM-xlyl-l}
    (A\tilde{x}^k + B\tilde{y}^k-b) -B(\tilde{y}^k-y^k) + (1/\beta) (\tilde{\lambda}^k-\lambda^k) = 0.
    \]
Combining \eqref{ADMM-xlyl-x}, \eqref{ADMM-xlyl-y} and \eqref{ADMM-xlyl-l}, and using the notation of \eqref{VI-FORM-2}, we obtain
\[\label{SC-PRSM-pre}
   \tilde{w}^k \in \Omega, \;\;
      \theta(u) -\theta(\tilde{u}^k)  + (w- \tilde{w}^k)^T  F(\tilde{w}^k)
          \ge  (v- \tilde{v}^k)^T Q(v^k-\tilde{v}^k), \;\; \forall w\in
          \Omega,
      \]
which corresponds to the prediction step (\ref{Frame-Q}) with $Q$ defined in (\ref{MatrixTQ}).

Moreover, note that $\lambda^{k+1}$ in  \eqref{S-ADMM-l2} can be represented as
\begin{eqnarray} \label{ADMM-SEl}
 \lambda^{k+1}
 & = & [{\lambda}^k -\mu(\lambda^k-\tilde{\lambda}^k)] -\mu \bigl[-\beta B(y^k-\tilde{y}^k) + \beta(Ax^{k+1}+By^k-b) \bigr] \nn \\
 & = & \lambda^k - \bigl[-\mu\beta B(y^k-\tilde{y}^k) + 2\mu(\lambda^k-\tilde{\lambda}^k) \bigr].
\end{eqnarray}
Thus, together with $y^{k+1}=\tilde{y}^k$, we have
$$
  \left(\begin{array}{c}
             y^{k+1} \\
              \lambda^{k+1} \end{array} \right)
             =
      \left(\begin{array}{c}
             y^{k} \\
              \lambda^{k} \end{array} \right)
                -    \left(\begin{array}{cc}
              I    &   0 \\
             -\mu\beta B   &   2\mu I_m \end{array} \right)  \left(\begin{array}{c}
             y^{k}-\tilde{y}^k \\
              \lambda^{k} - \tilde{\lambda}^k \end{array}
              \right).
             $$
This can be rewritten as the compact form
\[
     v^{k+1} = v^k -
      M (v^k-\tilde{v}^k),
     \]
which corresponds to the correction step (\ref{Frame-M}) with $M$ defined in (\ref{MatrixTQ}). \qquad {\normalsize$\Box$}

Hence, we can represent the SC-PRSM (\ref{S-ADMM-xlyl}) as the prediction-correction framework \eqref{Framework} with the corresponding matrices $Q$ and $M$ defined in (\ref{MatrixTQ}).

\subsubsection{Convergence conditions}

To verify that the matrices $Q$ and $M$ defined in (\ref{MatrixTQ}) satisfy the conditions (\ref{Frame-HG}) and hence ensure the convergence of the SC-PRSM (\ref{S-ADMM-xlyl}), we see that
$$  M^{-1} = \Biggl(\begin{array}{cc}
              I    &    0 \\
             \frac{1}{2}\beta B   &   \frac{1}{2\mu} I_m \end{array} \Biggr).
             $$
Hence, if we define
$$
  H = QM^{-1}=\Biggl(\!\!\!\begin{array}{cc}
                \beta B^TB& -\mu B^T\\
         -B  & \frac{1}{\beta} I_m \end{array}\!\!\! \Biggr)\!\Biggl(\!\!\!\begin{array}{cc}
              I    &    0 \\
             \frac{1}{2}\beta B   &   \frac{1}{2\mu} I_m \end{array} \!\!\!\Biggr)
         =\frac{1}{2} \Biggl(\!\!\!\begin{array}{cc}
              (2- \mu)\beta B^TB &    - B^T \\
            - B   &   \frac{1}{\mu\beta} I_m \end{array}\!\!\! \Biggr),
            $$
then $H$ is positive definite for any  $\mu\in(0,1)$ when $B$ is full column rank, and the
condition \eqref{Frame-H} is satisfied.

To check the condition \eqref{Frame-G}, i.e., the positiveness of $G:= Q^T+Q- M^THM $, we have
$$
   M^THM  = M^T Q  =  \Biggl(\!\!\begin{array}{cc}
              I    &    -\mu\beta B^T  \\
             0  &   2\mu I_m \end{array}\!\! \Biggr)\!\Biggl(\!\!\begin{array}{cc}
                \beta B^TB& -\mu B^T\\
         -B  & \frac{1}{\beta} I_m \end{array}\!\!\Biggr)  =  \Biggl(\!\! \begin{array}{cc}
              (1+\mu)\beta B^TB &   - 2\mu B^T \\
                 -2\mu B  & \frac{2\mu}{\beta} I_m
            \end{array}
           \!\!\Biggr).
     $$
Using \eqref{MatrixTQ} and the above equation, we have
$$ 
       G =Q^T+ Q- M ^THM   =  (1-\mu) \Biggl( \begin{array}{cc}
                \beta B^TB &   -  B^T \\
                 - B  & \frac{2}{\beta} I_m
            \end{array}
           \Biggr),
    $$
which implies that $G$ is positive definite for any $\mu\in (0,1)$ when $B$ is full column rank. Hence, the convergence conditions \eqref{Frame-HG} are satisfied for the strictly contractive PRSM (\ref{S-ADMM-xlyl}).

\section{Construction of new splitting contraction algorithms}  \label{Sec-3}

\setcounter{equation}{0}

As shown, the prediction-correction framework (\ref{Framework}), along with the generic convergence-guaranteeing conditions (\ref{Frame-HG}), represents a unified and abstract roadmap to convergence analysis for various splitting contraction algorithms in \cite{HeYuanCOA18,HeLWY,HeTaoYuanSIAMO,HeTaoYuanIMA,HeXuYuan}. In this section, in turn, we focus on how to specify the prediction-correction framework (\ref{Framework}) and the conditions (\ref{Frame-HG}). More specifically, with a given $Q$, we construct $H$ and $M$ to satisfy the conditions (\ref{Frame-HG}), and each pair of $H$ and $M$ can automatically generate a specific splitting contraction algorithm with provable convergence through the prediction-correction framework (\ref{Framework}). The constructive way is a principle of designing new application-tailored splitting contraction algorithms when concrete applications of the canonical convex programming model (\ref{Problem-LC}) are considered.

To design a new algorithm, it is rare to start from scratch; it is more often to start from a given coarse scheme that might be imperfect in theoretical or numerical aspects. Our discussion starts from the scenario where the matrix $Q$ is already determined by, e.g., a given coarse iterative scheme such as (\ref{DADMM}), (\ref{PD-Prediction}), or (\ref{DP-Prediction}), which will be delineated in the next two sections. With a given $Q$, our recept essentially only requires determining the norm matrix $H$ and the correction matrix $M$ for the prediction-correction framework (\ref{Framework}), while only the correction matrix $M$ needs to be specified to implement the resulting splitting contraction algorithm. Certainly, we prefer to construct such $M$ that can alleviate the resulting implementation/computation. We further assume that the given matrix $Q$ satisfies
$$   Q^T +  Q \succ  0. $$
Below, we give two specific principles of constructing the matrices $H$ and $M$ which can ensure the conditions (\ref{Frame-HG}).

\subsection{Construction from the condition $\boldsymbol{(\ref{Frame-H})}$}  \label{Sec3.1}

Note that the condition \eqref{Frame-H} can be rewritten as
  \[   \label{E-gHQD}    H= QM^{-1}.   \]
Since the norm matrix $H$ is required to be symmetric and positive definite, the condition (\ref{E-gHQD}) implies that $H$ should be representable in form of
  \[    \label{E-gHQDQ}     H= QD^{-1}Q^T, \]
in which the matrix $D$ is a undetermined positive definite matrix. Indeed, by comparing \eqref{E-gHQD} with \eqref{E-gHQDQ}, we know that
   $M^{-1}= D^{-1}Q^T$ and thus
  \[   \label{E-aDM}               M=Q^{-T} D.  \]
Hence, although the matrix $D$ in (\ref{E-gHQDQ}) is still unknown, choosing $M$ as (\ref{E-aDM}) can ensure the condition (\ref{Frame-H}).

Now, we investigate the restriction on $D$ to ensure the condition (\ref{Frame-G}) with the matrix $M$ given as (\ref{E-aDM}). Notice that
\[   \label{E5-MTHM}   M^T HM = \bigl(DQ^{-1}\bigr)\bigl(QD^{-1}Q^T\bigr) \bigl(Q^{-T}D\bigr)  = D.  \]
With (\ref{E5-MTHM}), then the condition (\ref{Frame-G}) is reduced to
\[ \label{E-Q-DG-2}  G:=Q^T+Q-M^THM=  Q^T+ Q-D \succ  0.\]
Hence, to ensure the condition (\ref{Frame-G}), the only restriction on the positive definite matrix $D$ in (\ref{E-gHQDQ}) is
 \[  \label{E-D}  0\prec D\prec Q^T+Q,  \]
In other words, whenever $Q$ is given and it satisfies  $   Q^T +  Q \succ  0$, then both $H$ and $M$ can be constructed via the following steps:
\[  \label{Chose-D}
    \left\{\begin{array}{rcl}
       HM &=& Q, \\
       M^THM&=& D.
       \end{array}\right. \quad \Longleftrightarrow  \quad
        \left\{\begin{array}{rcl}
       HM& =& Q,\\
       Q^TM&=& D.
       \end{array}\right.
       \quad \Longleftrightarrow  \quad
        \left\{\begin{array}{rcl}
       H &= &QD^{-1}Q^T, \\
       M &=& Q^{-T}D.
       \end{array}\right.  .  \]
Through this construction, both the conditions (\ref{Frame-M}) and (\ref{Frame-G}) are guaranteed to be satisfied. Note that once the matrix $D$ is chosen according to (\ref{E-D}), the matrices $H$, $M$ and $G$ are all uniquely determined. Then, with the specified matrix $M$ in (\ref{E-aDM}), the correction step (\ref{Frame-M}) and thus the prediction-correction framework (\ref{Framework}) is also specified as a concrete contraction splitting algorithm for the VI(\ref{VIu-O})-(\ref{Notation-uFO}).

\subsection{Construction from the condition $\boldsymbol{(\ref{Frame-G})}$}  \label{Sec3.2}

Alternatively, we can start from the condition \eqref{Frame-G} to construct the norm matrix $H$ and the correction matrix $M$. Again, with a given $Q$ satisfying $Q^T +  Q \succ 0$, we can choose the profit matrix $G$ such that
 \[  \label{E-G}        0  \prec  G \prec Q^T + Q.\]
Denote
   \[  \label{E-Delta}  \Delta = Q^T + Q  - G,\]
which is positive definite.  According to \eqref{Frame-G}, we know that the matrices $H$ and $M$ should satisfy
$$M^THM =\Delta.  $$
Recall the condition (\ref{Frame-H}): $HM=Q$. Thus, with a chosen $G$ satisfying (\ref{E-G}), $H$ and $M$ can be constructed vis the following steps:
 \[  \label{Chose-De}
    \left\{\begin{array}{rcl}
         M^THM&=& \Delta,\\
        HM &=& Q.
       \end{array}\right. \quad \Leftrightarrow  \quad
        \left\{\begin{array}{rcl}
            Q^TM&=& \Delta,\\
              HM& =& Q.
       \end{array}\right.
       \quad \Leftrightarrow  \quad
        \left\{\begin{array}{rcl}
            M &=& Q^{-T}\Delta,\\
              H &= &Q\Delta^{-1}Q^T.
       \end{array}\right.    \]
Then, with the constructed matrix $M$ in (\ref{Chose-De}), the correction step (\ref{Frame-M}) and thus the prediction-correction framework (\ref{Framework}) can also be specified as a concrete splitting contraction algorithm for the VI(\ref{VIu-O})-(\ref{Notation-uFO}).
Again, with a given $G$ satisfying (\ref{E-G}), the matrices $H$ and $M$ are both uniquely determined.


%

\subsection{Remarks}\label{Sec3.3}

It is interesting to observe that the proposed two construction strategies can be related via the relationship
 \[\label{DG-QTQ}   D\succ 0, \quad G\succ 0, \quad \hbox{and}\quad   D + G = Q^T+Q. \]
Hence, once $D$ is chosen for the construction strategy in Section \ref{Sec3.1}, the corresponding $G$ given by (\ref{DG-QTQ}) can be used for the construction strategy in Section \ref{Sec3.2}, and vice versa.

Our discussions above show that, with a given prediction matrix $Q$ satisfying $ Q^T +  Q \succ  0$, once $D$ (Resp., $G$) is chosen according to (\ref{E-D}) (Resp., (\ref{E-G})), the matrices $H$ and $M$ can be determined as analyzed and thus a concrete contraction splitting algorithm can be specified for the convex programming problem (\ref{Problem-LC}) via the prediction-correction framework (\ref{Framework}). Technically, there are infinitely many such choices subject to (\ref{DG-QTQ}). For example, we can choose
$$
D= \alpha[Q^T+Q]  \quad \hbox{and}\quad G= (1-\alpha)[Q^T + Q], \;\; \;\;\alpha\in(0,1).
 $$
We will elaborate on the choice $D=  G= \frac{1}{2}[ Q^T+Q]$ in Section \ref{Sec4.3}.

\subsection{Implementation of the correction step (\ref{Frame-M})}

Note that the correction step \eqref{Frame-M} can be rewritten as
 $$   Q^T(v^{k+1} - v^k)  = Q^T      M (v^k-\tilde{v}^k). $$
To implement the correction step (\ref{Frame-M}) with the constructed two choices for $M$, i.e., $M=Q^{-T}D$ in (\ref{E-aDM})
and  $M=Q^{-T} \Delta$ in (\ref{Chose-De}), we need to solve one of the following systems of equations:
\[\label{QD}     Q^T (v^{k+1} - v^k) = D(\tilde{v}^k - v^k),    \]
 and
 \[\label{QG}
 Q^T (v^{k+1} - v^k) = \Delta(\tilde{v}^k - v^k).    \]
Hence, although $D$ and $G$ (thus $\Delta$) can be chosen arbitrarily with the only constraint (\ref{E-D}) or (\ref{E-G}), it is preferred to choose some model-tailored ones that can favor solving the systems of equations (\ref{QD}) or (\ref{QG}) more efficiently. Irrational choices that make the correction step (\ref{Frame-M}) complicated should be generally avoided. Some examples will be discussed in the next sections when specific cases of the canonical convex programming model (\ref{Problem-LC}) are considered.

\section{Application to three-block separable convex optimization}  \label{Section-S}

\setcounter{equation}{0}

In this section, we apply the strategies proposed in Sections \ref{Sec3.1}, \ref{Sec3.2} and \ref{Sec3.3} to a separable convex optimization problem, and showcase how to construct the norm matrix $H$ and the correction matrix $M$ when the matrix $Q$ is given.

\subsection{Model}

We consider the three-block separable convex optimization model with linear constraints
\[  \label{Problem-A3}
    \min \{\theta_1(x) + \theta_2(y) +\theta_3(z) \, | \,  Ax + By +Cz= b, \, x\in {\cal X}, y\in{\cal Y}, z\in{\cal Z} \},
   \]
where $\theta_i: {\Re}^{n_i}\to {\Re} \;(i=1,2,3)$ are closed, proper, and convex but not necessarily smooth functions; ${\cal X} \subseteq \Re^{n_1} $, ${\cal Y} \subseteq \Re^{n_2}$ and ${\cal Z} \subseteq \Re^{n_3}$ are closed convex sets; $A \in \Re^{m\times n_1}$, $B \in \Re^{m\times n_2}$ and  $C \in \Re^{m\times n_3}$; $b\in \Re^m$; and $n_1+n_2 + n_3 = n$. Clearly, it is a special case of the canonical convex programming problem (\ref{Problem-LC}), and the VI (\ref{VIu-O})-(\ref{Notation-uFO}) can be specified as the following:
\[  \label{VI-AS2}
      w^*\in \Omega, \quad \theta(u) -\theta(u^*) + (w-w^*)^T F(w^*) \ge 0, \quad \forall \,  w\in
      \Omega,
     \]
where
\begin{subequations} \label{VI-FORM-A2}
\[
      w = \left(\begin{array}{c}
                     x\\  y\\  z\\  \lambda \end{array} \right), \quad
                     u= \left(\begin{array}{c}
                     x\\  y\\z  \end{array} \right), \quad  F(w) = \left(\begin{array}{c}
            - A^T\lambda \\
             -B^T\lambda \\
             -C^T\lambda \\
             Ax + By +Cz -b \end{array} \right),
   \]
with
\[
   \theta(u)=\theta_1(x) +
                     \theta_2(y)  +\theta_3(z),  \qquad
            \Omega={\cal X} \times {\cal Y} \times {\cal Z} \times \Re^m.
   \]
\end{subequations}

Let the augmented Lagrangian function of the model \eqref{Problem-A3}  be
\[\label{ALF-3} {\cal L}_{\beta}^{[3]}(x,y,z,\lambda) =\theta_1(x) + \theta_2(y) + \theta_3(z) - \lambda^T(Ax
 +By +Cz -b) +\frac{\beta}{2}\|Ax
 +By +Cz -b\|^2
   \]
with $\lambda \in \Re^m$ the Lagrange multiplier and $\beta > 0$ the penalty parameter. With the success of the ADMM (see, e.g.\cite{FG,Glow84,Glowinski89}), it is natural to consider directly extending the ADMM (\ref{ADMM-2}) and splitting the augmented Lagrangian function ${\cal L}_{\beta}^{[3]}(x,y,z,\lambda)$ in (\ref{ALF-3}) three times by the Gauss-Seidel manner. That is, consider the scheme
\[\label{DADMM}
  \hskip-0.09cm\left\{\begin{array}{l}
    x^{k+1} \in \arg\min\bigl\{  {\cal L}_{\beta}^{[3]}(x,y^k,z^k,\lambda^k) \;|\;  x\in {\cal X}\bigr\},  \\[0.2cm]
    {y}^{k+1} \in \arg\min\bigl\{ {\cal L}_{\beta}^{[3]}({x}^{k+1}, y,z^k,\lambda^k)
                          \;|\;  y\in {\cal Y}\bigr\},  \\[0.2cm]
    {z}^{k+1} \in \arg\min\bigl\{  {\cal L}_{\beta}^{[3]}({x}^{k+1},{y}^{k+1},z,\lambda^k)
                          \;|\;  z\in {\cal Z}\bigr\}, \\[0.2cm]
    {\lambda}^{k+1}= \lambda^k -  \bigl(A{x}^{k+1} +B{y}^{k+1} + Cz^{k+1} -b\bigr).
  \end{array} \right.
  \]
However, the splitting scheme (\ref{DADMM}) is coarse in sense of that its convergence is not guaranteed as shown in \cite{CHYY}. Thus, to render convergence, either the scheme (\ref{DADMM}) should be appropriately adjusted or stronger conditions on functions/coefficient matrices/penalty parameters should be additionally assumed. To develop algorithms at the root of the ADMM (\ref{ADMM-2}), we can consider correcting the output of (\ref{DADMM}) by certain correction steps as what we did in \cite{HeYuanCOA18,HeTaoYuanSIAMO,HeTaoYuanIMA,HeXuYuan}. Below we show how to apply the construction strategies in Section \ref{Sec-3} to generate splitting contraction algorithms by taking advantage of the coarse splitting scheme (\ref{DADMM}).

\subsection{Discerning the prediction matrix $Q$}

Our construction starts from the coarse splitting scheme (\ref{DADMM}) which can be rewritten as the prediction step (\ref{Frame-Q}) and hence the corresponding prediction matrix $Q$ can be discerned. For this purpose, we first consider the subproblems related to the primal variables in (\ref{DADMM}), and rewrite them as
 $\tilde{u}^k
   =(\tilde{x}^k,\tilde{y}^k,\tilde{z}^k)\in {\cal X}\times{\cal Y}\times{\cal Z}$. Namely, we have
\[\label{DADMM-tilde}
  \hskip-0.09cm\left\{\begin{array}{l}
    \tilde{x}^k \in \arg\min\bigl\{  {\cal L}_{\beta}^{[3]}(x,y^k,z^k,\lambda^k) \;|\;  x\in {\cal X}\bigr\},  \\[0.2cm]
      \tilde{y}^k \in \arg\min\bigl\{ {\cal L}_{\beta}^{[3]}(\tilde{x}^k, y,z^k,\lambda^k)
                          \;|\;  y\in {\cal Y}\bigr\},  \\[0.2cm]
      \tilde{z}^k \in \arg\min\bigl\{  {\cal L}_{\beta}^{[3]}(\tilde{x}^k,\tilde{y}^k,z,\lambda^k)
                          \;|\;  z\in {\cal Z}\bigr\}.
  \end{array} \right.
  \]
Ignoring some constant terms, we can rewrite the formula above as
  \[ \label{PD-Preu}
  \hskip-0.09cm\left\{\begin{array}{l}
    \tilde{x}^k \in \arg\min\bigl\{   \theta_1(x) -  x^TA^T\lambda^k
                         +  \frac{\beta}{2} \|Ax+By^k+Cz^k-b\|^2
                          \;|\;  x\in {\cal X}\bigr\},  \\[0.2cm]
      \tilde{y}^k \in \arg\min\bigl\{   \theta_2(y) -  y^TB^T\lambda^k
                         +  \frac{\beta}{2} \|A\tilde{x}^k+  By  +Cz^k-b\|^2
                          \;|\;  y\in {\cal Y}\bigr\},  \\[0.2cm]
      \tilde{z}^k \in \arg\min\bigl\{   \theta_3(z) -  z^TC^T\lambda^k
                         +  \frac{\beta}{2} \|A\tilde{x}^k+B\tilde{y}^k+Cz-b\|^2
                          \;|\;  z\in {\cal Z}\bigr\}.
  \end{array} \right.
  \]
Then, according to Lemma \ref{CP-TF}, we have $\tilde{u}^k\in {\cal U}$ and
 \[
\left\{ \begin{array}{l}         \theta_1(x) - \theta_1(\tilde{x}^k)  + (x - \tilde{x}^k)^T\bigl\{-A^T{\lambda}^k \\
  \hskip 4.5cm +\beta A^T\! \bigl(A\tilde{x}^k+By^k + Cz^k-b\bigr)\bigr\} \ge 0,  \;\; \forall \, x\in {\cal X},\\[0.1cm]
       \theta_2(y) - \theta_2(\tilde{y}^k)  + (y - \tilde{y}^k)^T\bigl\{-B^T{\lambda}^k \\
  \hskip 4.5cm +\beta   B^T\! \bigl(A\tilde{x}^k+B\tilde{y}^k + Cz^k-b\bigr)\bigr\} \ge 0,  \;\;  \forall \,  y\in {\cal Y}, \\[0.1cm]
     \theta_3(z) - \theta_3(\tilde{z}^k)  + (z - \tilde{z}^k)^T\bigl\{-C^T{\lambda}^k \\
  \hskip 4.5cm +\beta  C^T\! \bigl(A\tilde{x}^k+B\tilde{y}^k + C\tilde{z}^k-b\bigr)\bigr\} \ge 0,  \;\;  \forall \,  z\in {\cal Z}.
                \end{array}\right.
   \]
Defining
\[ \label{PD-Prel}   \tilde{\lambda}^k = \lambda^k  - \beta\bigl(A\tilde{x}^k+By^k + Cz^k-b\bigr),\]
and using the VI form \eqref{VI-FORM-A2}, we have  $\tilde{w}^k\in \Omega$ and
  \[ \label{DP-C}
      \left\{ \begin{array}{l}
       \theta_1(x) - \theta_1(\tilde{x}^k)  + (x - \tilde{x}^k)^T\bigl\{ \underline{-A^T\tilde{\lambda}^k}\bigr\} \ge 0,  \quad \forall \, x\in {\cal X},\\[0.1cm]
        \theta_2(y) - \theta_2(\tilde{y}^k)  + (y - \tilde{y}^k)^T
        \bigl\{ \underline{-B^T\tilde{\lambda}^k}    +     \beta  B^T\!B(\tilde{y}^k-y^k) \bigr\} \ge 0, \quad \forall \, y\in {\cal Y}, \\[0.1cm]
     \theta_3(z) - \theta_3(\tilde{z}^k)  + (z - \tilde{z}^k)^T
     \Biggl\{\!\begin{array}{c}{\underline{-C^T\tilde{\lambda}^k}}
           + \beta C^T\!B(\tilde{y}^k-y^k)\\
           \beta   C^T\!C(\tilde{z}^k-z^k)  \end{array}\! \Biggr\} \ge 0, \quad \forall \, z\in {\cal Z}, \\[0.1cm]
   (\lambda - \tilde{\lambda}^k)^T\Biggl\{\begin{array}{l}
                 ({\underline{A\tilde{x}^k + B\tilde{y}^k  + C\tilde{z}^k-b}})  \\
                       - \!B(\tilde{y}^k-y^k) \! -\!C(\tilde{z}^k-z^k) + \dfrac{1}{\beta}  (\tilde{\lambda}^k-\lambda^k)
                         \end{array}\!\Biggr\}\ge 0, \quad  \forall \, \lambda\in {\Lambda}.
  \end{array} \right.
              \]
The sum of the underline parts of \eqref{DP-C} is exactly $F(\tilde{w}^k)$, where $F(\cdot)$ is defined in (\ref{VI-FORM-A2}). Thus, we have
\[  \label{DP-PreD-F }
  \tilde{w}^k\in \Omega, \;\;\;   \theta(u) - \theta(\tilde{u}^k) +  (w- \tilde{w}^k )^T
       F(\tilde{w}^k) \ge   (v- \tilde{v}^k )^T  Q(  v^k -\tilde{v}^k),  \quad \forall \, w \in\Omega,
\]
where the prediction matrix is
     \[ \label{DP-Q}
      Q=  \left(\begin{array}{ccc}
          \beta  B^TB  &   0                  &        0      \\
          \beta    C^TB  & \beta  C^TC &       0         \\
           -B     &   -C    &  \dfrac{1}{\beta}I_m
   \end{array} \right).
    \]

Moreover, for the prediction matrix $Q$ in (\ref{DP-Q}) which is determined by the coarse splitting scheme (\ref{DADMM}), we have
\[\label{QTQ}   Q^T +Q =
        \left(\begin{array}{ccc}
            2\beta   B^T\!B&   \beta  B^T\!C    &       -B^T\\
             \beta   C^T\!B   &   2\beta  C^T\!C      &  -C^T \\
               -B   &     -C  &  \frac{2}{\beta} I_m
   \end{array} \right),  \]
which is positive definite whenever $B$ and $C$ are full column rank.

\subsection{Constructing the correction matrix $M$}

With the prediction matrix $Q$ given in (\ref{DP-Q}), the prediction-correction framework (\ref{Framework}) can be specified as a concrete algorithm for the model (\ref{Problem-A3}) once the correction step (\ref{Frame-M}) is specified. Now, we showcase how to specify the correction step (\ref{Frame-M}) by the construction strategies discussed in Sections \ref{Sec3.1}, \ref{Sec3.2}  and \ref{Sec3.3}. Note that $v=(y,z,\lambda)$ below.

\subsubsection{Construction 1}

Based on (\ref{DP-Q}) and (\ref{QTQ}), and following the strategy in Section \ref{Sec3.1}, we can choose
   \[  \label{Example-D}   D=  \left(\begin{array}{ccc}
              \nu  \beta  B^T\!B&   0    &       0 \\
                   0 &    \nu\beta C^T\!C     &  0 \\
                       0   &    0  &\dfrac{1}{\beta }  I_m
   \end{array} \right)   \]
with $0<\nu<1$, which is positive definite whenever $B$ and $C$ are both full column rank. Recall the correction matrix $M$ in (\ref{E-aDM}). Then, a concrete splitting contraction algorithm for (\ref{Problem-A3}) can be generated as below.
\begin{center}\fbox{\begin{minipage}{15.5cm}

\begin{subequations}
\noindent{\bf \normalsize Algorithm 1 for the model (\ref{Problem-A3})}

[Prediction Step.] Obtain $(\tilde{x}^k,\tilde{y}^k,\tilde{z}^k)$ via the direct extension of the ADMM (\ref{DADMM-tilde})   and \par
   \qquad \qquad \qquad \qquad  define  $\tilde{\lambda}^k$ by (\ref{PD-Prel}).

\smallskip

[Correction Step.] $Q^T (v^{k+1} -v^k) = D (\tilde{v}^k -v^k)$.
\end{subequations}

\end{minipage}
}
\end{center}


For the correction step $Q^T (v^{k+1} -v^k) = D (\tilde{v}^k -v^k)$, we know that
  $$  Q^T=  \left(\begin{array}{ccc}
              \! \beta B^T\!B\!& \! \beta B^T\!C\!   &    -B^T \\
                   0  &\! \beta C^T\!C\!     &  -C^T \\
                          0   &    0   &\!\frac{1}{\beta}  I_m   \!
   \end{array} \right)  = \left(\begin{array}{ccc}
                \beta B^T&    0    &    0 \\
                   0  & \beta  C^T      &  0 \\
                          0   &    0   & \frac{1}{\beta}    I_m
   \end{array} \right) {\left(\begin{array}{ccc}
              \!  B\!& \! C\!   &    -\frac{1}{\beta}I_m\\
                   0  &  \! C\!     &  -\frac{1}{\beta} I_m \\
                     0   &  0  & \! I_m\!
   \end{array} \right)}, $$
and
 $$  D=  \left(\begin{array}{ccc}
                \!\nu\beta  B^T\!B\!&  0   &    0 \\
                   0  &\! \nu \beta  C^T\!C\!     &  0 \\
                   0   &    0   &\! \frac{1}{\beta} I_m   \!
   \end{array} \right)  =\left(\begin{array}{ccc}
                \beta B^T&    0    &    0 \\
                   0  & \beta  C^T      &  0 \\
                          0   &    0   & \frac{1}{\beta}    I_m
   \end{array} \right){ \left(\begin{array}{ccc}
              \!  \nu B\!&   0   &    0\\
                   0  &  \!\nu C\!     &  0 \\
                   0  &    0    & \! I_m\!
   \end{array} \right)}. $$
That is, $Q^T$ and $D$ have a common matrix in their factorization forms above. Hence, to implement the correction step (\ref{Frame-M}), i.e., $Q^T (v^{k+1} -v^k) = D (\tilde{v}^k -v^k)$, essentially we only need to consider the even easier equation
     \[  \label{QTD-Av}   \left(\begin{array}{ccc}
              \!  B\!& \! C\!   &    -\frac{1}{\beta}I_m\\
                   0  &  \! C\!     &  -\frac{1}{\beta} I_m \\
                     0   &  0  & \! I_m\!
   \end{array} \right)
       \left(\begin{array}{ccc}
            y^{k+1}-y^k  \\
             z^{k+1} - z^k \\
              \lambda^{k+1}-\lambda^k
   \end{array} \right) =  \left(\begin{array}{ccc}
               \nu B&    0    &    0 \\
                   0  &   \nu C      &  0 \\
                          0   &    0   &  I_m
   \end{array} \right)
       \left(\begin{array}{ccc}
            \tilde{y}^{k}-y^k  \\
             \tilde{z}^{k} - z^k \\
              \tilde{\lambda}^{k}-\lambda^k
   \end{array} \right). \]
Moreover, iterations of the specified Algorithm 1 can be executed in terms of $(By^{k},Cz^{k}, \lambda^{k})$ because it is sufficient to keep $By^{k}$ and $Cz^{k}$, rather than $y^k$ and $z^k$, to execute the $(k+1)$-th iteration. The variables $y^k$ and $z^k$ need to be solved once only at the last iteration. Hence, with the choice of $D$ in (\ref{Example-D}), implementing the resulting correction step (\ref{Frame-M}) essentially only requires solving the equation (\ref{QTD-Av}) in terms of $(By^{k},Cz^{k}, \lambda^{k})$ , which is extremely easy.

\subsubsection{Construction 2}

Based on (\ref{DP-Q}) and (\ref{QTQ}), and following the strategy in Section \ref{Sec3.2}, we can choose
  \[ \label{Example-G}  G=  \left(\begin{array}{ccc}
              \nu  \beta B^T\!B&   0    &       0 \\
                   0 &    \nu\beta  C^TC     &  0 \\
                       0   &    0  &\frac{1}{\beta }  I_m \end{array} \right), \]
with $\nu\in(0,1)$, which can be guaranteed to be positive definite whenever $B$ and $C$ are full column rank. Note that the matrix $G$ in (\ref{Example-G}) is precisely the matrix $D$ defined in \eqref{Example-D}. Furthermore, we have
 $$ \Delta =Q^T+Q- G =  \left(\begin{array}{ccc}
             (2-\nu)\beta B^T\!B  &  \beta   B^TC    &     -B^T \\
               \beta C^TB   &   (2-\nu)\beta  C^TC   &   -C^T \\
                          -B   &   -C    & \frac{1}{\beta } I_m
   \end{array} \right). $$
Recall the correction matrix $M$ in (\ref{Chose-De}). Then, another contraction splitting algorithm for (\ref{Problem-A3}) can be generated as below.

\begin{center}\fbox{\begin{minipage}{15.5cm}

\begin{subequations}
\noindent{\bf \normalsize Algorithm 2 for the model (\ref{Problem-A3})}

[Prediction Step.] Obtain $(\tilde{x}^k,\tilde{y}^k,\tilde{z}^k)$ via the direct extension of the ADMM (\ref{DADMM-tilde}) and \par
    \qquad \qquad \qquad \qquad  define  $\tilde{\lambda}^k$ by (\ref{PD-Prel}).
\smallskip

[Correction Step.] $Q^T (v^{k+1} -v^k) = \Delta  (\tilde{v}^k -v^k)$.
\end{subequations}

\end{minipage}
}
\end{center}

For the correction step $Q^T (v^{k+1} -v^k) = \Delta  (\tilde{v}^k -v^k)$, we know that
  $$   Q^T=  \left(\begin{array}{ccc}
              \! \beta B^T\!B\!& \! \beta B^T\!C\!   &    -B^T \\
                   0  &\! \beta C^T\!C\!     &  -C^T \\
                          0   &    0   &\!\frac{1}{\beta}  I_m   \!
   \end{array} \right)  =\left(\begin{array}{ccc}
                \beta B^T&    0    &    0 \\
                   0  & \beta  C^T      &  0 \\
                          0   &    0   & \frac{1}{\beta}    I_m
   \end{array} \right) {\left(\begin{array}{ccc}
              \!  B\!& \! C\!   &    -\frac{1}{\beta}I_m\\
                   0  &  \! C\!     &  -\frac{1}{\beta} I_m \\
                     0   &  0  & \! I_m\!
   \end{array} \right)}, $$
and
 \begin{eqnarray*}
   \Delta
   & =&   \left(\begin{array}{ccc}
             (2-\nu)\beta B^T\!B  &  \beta   B^TC    &     -B^T \\
               \beta C^TB   &   (2-\nu)\beta  C^TC   &   -C^T \\
                          -B   &   -C    & \frac{1}{\beta } I_m
   \end{array} \right)    \nn \\
    &= &   \left(\begin{array}{ccc}
                \beta B^T&    0    &    0 \\
                   0  & \beta  C^T      &  0 \\
                          0   &    0   & \frac{1}{\beta}    I_m
   \end{array} \right){
     \left(\begin{array}{ccc}
             (2-\nu)B &    C             &  -\frac{1}{\beta} I_m \\
                B           &   (2-\nu) C  &  - \frac{1}{\beta}I_m \\
                 -\beta B      &     -\beta C              &  I_m
   \end{array} \right)}.
   \end{eqnarray*}
That is, $Q^T$ and $\Delta$ have a common matrix in their factorization forms above. Hence, to implement the correction step (\ref{Frame-M}), i.e., $Q^T (v^{k+1} -v^k) = \Delta (\tilde{v}^k -v^k)$, essentially we only need to consider the even easier equation
\[   \label{QTDe-Av}
   \left(\!\begin{array}{ccc}
              \!  B\!& \! C\!   &    -\frac{1}{\beta}I_m\\
                   0  &  \! C\!     &  -\frac{1}{\beta} I_m \\
                     0   &  0  & \! I_m\!
   \end{array}\! \right)
       \left(\begin{array}{ccc}
            y^{k+1}-y^k  \\
             z^{k+1} - z^k \\
              \lambda^{k+1}-\lambda^k
   \end{array} \right)
   =
     \left(\!\!\begin{array}{ccc}
             (2-\nu)B &    C             &  -\frac{1}{\beta} I_m \\
                B           &   (2-\nu) C  &  - \frac{1}{\beta}I_m \\
                 -\beta B      &     -  \beta C      &  I_m
   \end{array}\!\! \right)    \left(\begin{array}{ccc}
            \tilde{y}^{k}-y^k  \\
             \tilde{z}^{k} - z^k \\
              \tilde{\lambda}^{k}-\lambda^k
   \end{array} \right).
   \]
Similar as (\ref{QTD-Av}), with the choice of $G$ in (\ref{Example-G}), implementing the resulting correction step (\ref{Frame-M}) essentially only requires solving the equation (\ref{QTDe-Av}) in terms of $(By^{k},Cz^{k}, \lambda^{k})$, which is extremely easy.

\subsubsection{Construction 3}\label{Sec4.3}

Recall the relationship between the matrices $D$ and $G$ in (\ref{DG-QTQ}), and $Q^T+Q$ given in (\ref{QTQ}). Essentially, the proposed construction strategies in Sections \ref{Sec3.1} and \ref{Sec3.2} take the same matrix
$$ \left(\begin{array}{ccc}
              \nu  \beta B^T\!B&   0    &       0 \\
                   0 &    \nu\beta  C^TC     &  0 \\
                       0   &    0  &\frac{1}{\beta }  I_m \end{array} \right) $$
as $D$ and $G$, respectively, and then the other one is determined by (\ref{DG-QTQ}). As mentioned in Section \ref{Sec3.3}, any other choice of $D$ and $G$ subject to the relationship (\ref{DG-QTQ}) is also eligible. Let us consider the following specific one:
   \[ \label{Example-DG}    D=G =\frac{1}{2}\bigl[Q^T +Q \bigr] = \left(\begin{array}{ccc}
            \beta   B^T\!B& \frac{1}{2}  \beta  B^T\!C      &   - \frac{1}{2}B^T\\[0.1cm]
              \frac{1}{2}\beta   C^T\!B   &   \beta  C^T\!C &  - \frac{1}{2}C^T \\[0.1cm]
               - \frac{1}{2}B   &      \frac{1}{2}-C  &  \frac{}{\beta} I_m
   \end{array} \right) , \]
which are both positive definite whenever $B$ and $C$ are full column rank. Recall the correction matrix $M$ in (\ref{Chose-De}). Then, one more contraction splitting algorithm for (\ref{Problem-A3}) can be generated as below.

\begin{center}\fbox{\begin{minipage}{15.5cm}

\begin{subequations}
\noindent{\bf \normalsize Algorithm 3 for the model (\ref{Problem-A3})}

[Prediction Step.] Obtain $(\tilde{x}^k,\tilde{y}^k,\tilde{z}^k)$ via the direct extension of the ADMM (\ref{DADMM-tilde}) and \par
    \qquad \qquad \qquad \qquad  define  $\tilde{\lambda}^k$ by (\ref{PD-Prel}).
\smallskip

[Correction Step.] $Q^T (v^{k+1} -v^k) = \frac{1}{2}[Q^T +Q]  (\tilde{v}^k -v^k)$.
\end{subequations}

\end{minipage}
}
\end{center}

For the correction step $Q^T (v^{k+1} -v^k) = \frac{1}{2}[Q^T +Q]  (\tilde{v}^k -v^k)$, we know that
  $$   Q^T=  \left(\begin{array}{ccc}
              \! \beta B^T\!B\!& \! \beta B^T\!C\!   &    -B^T \\
                   0  &\! \beta C^T\!C\!     &  -C^T \\
                          0   &    0   &\!\frac{1}{\beta}  I_m   \!
   \end{array} \right)  = \left(\begin{array}{ccc}
                \beta B^T&    0    &    0 \\
                   0  & \beta  C^T      &  0 \\
                          0   &    0   & \frac{1}{\beta}    I_m
   \end{array} \right) {\left(\begin{array}{ccc}
              \!  B\!& \! C\!   &    -\frac{1}{\beta}I_m\\
                   0  &  \! C\!     &  -\frac{1}{\beta} I_m \\
                     0   &  0  & \! I_m\!
   \end{array} \right)}, $$
and

$$
 \frac{1}{2}[Q^T+Q]
   =  \left(\!\!\begin{array}{ccc}
            \beta   B^T\!B& \frac{1}{2}  \beta  B^T\!C      &   - \frac{1}{2}B^T\\[0.1cm]
              \frac{1}{2}\beta   C^T\!B   &   \beta  C^T\!C &  - \frac{1}{2}C^T \\[0.1cm]
               - \frac{1}{2}B   &      \frac{1}{2}-C  &  \frac{1}{\beta} I_m
   \end{array}\! \right)  =   \left(\!\!\begin{array}{ccc}
                \beta B^T&    0    &    0 \\
                   0  & \beta  C^T      &  0 \\
                          0   &    0   & \frac{1}{\beta}    I_m
   \end{array}\!\!\right)\!\!
     \left(\!\!\begin{array}{ccc}
             B &   \frac{1}{2} C             &  -\frac{1}{2\beta} I_m \\
              \frac{1}{2}   B           &   C  &  - \frac{1}{2\beta}I_m \\
                 -\frac{1}{2}\beta B      &     -\frac{1}{2}\beta C              &  I_m
   \end{array}\!\! \right).  $$
That is, $Q^T$ and $\frac{1}{2}[Q^T +Q]$ have a common matrix in their factorization forms above. Hence, to implement the correction step (\ref{Frame-M}), i.e., $Q^T (v^{k+1} -v^k) = \frac{1}{2}[Q^T +Q]  (\tilde{v}^k -v^k)$, essentially we only need to consider the even easier equation

\[  \label{QTDG-Av}
    \left(\begin{array}{ccc}
              \!  B\!& \! C\!   &    -\frac{1}{\beta}I_m\\
                   0  &  \! C\!     &  -\frac{1}{\beta} I_m \\
                     0   &  0  & \! I_m\!
   \end{array} \right)
       \left(\begin{array}{ccc}
            y^{k+1}-y^k  \\
             z^{k+1} - z^k \\
              \lambda^{k+1}-\lambda^k
   \end{array} \right)
   =
     \left(\begin{array}{ccc}
             B &   \frac{1}{2} C             &  -\frac{1}{2\beta} I_m \\
              \frac{1}{2}   B           &   C  &  - \frac{1}{2\beta}I_m \\
                 -\frac{1}{2}\beta B      &     -\frac{1}{2}\beta C              &  I_m
   \end{array} \right) \left(\begin{array}{ccc}
            \tilde{y}^{k}-y^k  \\
             \tilde{z}^{k} - z^k \\
              \tilde{\lambda}^{k}-\lambda^k
   \end{array} \right).\]
   Similar as (\ref{QTD-Av}) and (\ref{QTDe-Av}), implementing the resulting correction step (\ref{Frame-M}) essentially only requires solving the equation (\ref{QTDG-Av}) in terms of $(By^{k},Cz^{k}, \lambda^{k})$, which is extremely easy.

\subsubsection{Summary}

We have discussed three concrete strategies for constructing ADMM-based splitting contraction algorithms for the model (\ref{Problem-A3}), with the same prediction step determined by the coarse splitting scheme (\ref{DADMM-tilde}) with (\ref{PD-Prel}). These three strategies differ in how to choose the matrices $D$ and $G$ subject to (\ref{DG-QTQ}); accordingly specifications of the correction step (\ref{Frame-M}) are different. It is easy to verify that Algorithm 1 corresponds to the algorithm proposed in \cite{HeTaoYuanSIAMO}. Certainly, any other choice of $D$ and $G$ in accordance with (\ref{DG-QTQ}) leads to another splitting contraction algorithm for the model (\ref{Problem-A3}) whose prediction step remains unchanged as that in Algorithms 1-3, while there are infinitely many such choices.

\section{Application to multiple block separable convex optimization}  \label{SecM-5}

\setcounter{equation}{0}

In this section, we extend the analysis in Section \ref{Section-S} to more general and complicated convex programming problems. We recall a recent work of ours \cite{HeXuYuan}, and will show that more new algorithms with similar advantages as those in \cite{HeXuYuan} can be presented by following the proposed constructions strategies in Section \ref{Sec-3}.

\subsection{Model}

We consider the generic multiple block convex programming problem with both linear equality and inequality constraints:

\[  \label{Problem-m}
  \min \Bigl\{ \sum_{i=1}^{p} \theta_i(x_i)   \;\big|\;   \sum_{i=1}^{p} A_ix_i=b\ (\hbox{or} \ge b) ,  \;\;  x_i\in {\cal X}_i \Bigr\},\]
where $\theta_i: {\Re}^{n_i}\to {\Re}\; (i=1,\ldots, p)$ are closed, proper, and convex but not necessarily smooth functions; ${\cal
X}_i\subseteq \Re^{n_i}\; (i=1,\ldots, p)$ are closed convex sets;
$A_i\in \Re^{m\times n_i}\; (i=1,\ldots, p)$ are given matrices; $b\in \Re^m$; $\sum_{i=1}^{p}n_i=n$; and $p\ge 3$.
Note that the more general model (\ref{Problem-m}) differs from the special three-block separable convex programming problem (\ref{Problem-A3}) in that $p\ge3$ and that both linear equality and inequalities are considered.

Let the Lagrangian function of \eqref{Problem-m} be
\[ \label{Lagrange-F}
  L(x_1,\ldots,x_p,\lambda) =  \sum_{i=1}^{p}\theta_i(x_i) -\lambda^T\Bigl(\sum_{i=1}^{p} A_ix_i-b\Bigr),\]
with $\lambda$ the Lagrange multiplier. We know that $\lambda\in \Re^m$ or $\lambda\in \Re^m_+$ if $\sum_{i=1}^{p} A_ix_i=b $ or $\sum_{i=1}^{p} A_ix_i\ge b$ is considered in the model (\ref{Problem-m}), respectively. The optimality condition of (\ref{Problem-m}) can be written as the VI (\ref{VIu-O})-(\ref{Notation-uFO}) with the following specifications:
\begin{subequations} \label{VI-pFORM}
\[  \label{VI-pFORM-Q}
    w^*\in\Omega, \quad \theta(x) - \theta(x^*) + (w-w^*)^T F(w^*) \ge0, \quad
\forall\, w\in \Omega,
    \]
where
\[  \label{VI-pFORM-F}    w=\left(\!\!\begin{array}{c}
             x_1\\
         \vdots \\
            x_{p} \\[0.1cm]
           \lambda
             \end{array}\!\! \right), \quad
   x=\left(\!\!\begin{array}{c}
             x_1\\
         \vdots \\
            x_{p} \\
             \end{array}\!\! \right),
    \quad   \theta(x) = \sum_{i=1}^{p} \theta_i(x_i), \quad
       F(w) = \left(\!\!\begin{array}{c}
    - A_1^T\lambda \\
       \vdots \\
    -A_{p}^T\lambda \\[0.1cm]
       \sum_{i=1}^{p} A_ix_i-b
    \end{array}\!\! \right),
  \]
  and
  $$      \Omega = \prod_{i=1}^p {\cal X}_i \times \Lambda \quad \hbox{with}\quad   \Lambda  =\left\{ \begin{array}{ll}
               \Re^m,           &   \hbox{if   $\sum_{i=1}^{p} A_ix_i =  b$} ,   \\[0.2cm]
                \Re^m_+,    &        \hbox{if   $\sum_{i=1}^{p} A_ix_i\ge b$}.
                \end{array}        \right.  $$
\end{subequations}

\subsection{Discerning the prediction matrix $Q$}

To design an algorithm for (\ref{Problem-m}), it is natural to consider extending the ADMM (\ref{ADMM-2}) and splitting the augmented Lagrangian function $L(x_1,\ldots,x_p,\lambda)$ in (\ref{Lagrange-F}) to obtain easier subproblems. In \cite{HeXuYuan}, two different ways were suggested and the resulting subproblems were used as prediction steps for the prediction-correction framework (\ref{Framework}). We recall these results; accordingly the prediction matrix $Q$ becomes clear.

\subsubsection{Primal-dual order}

We first consider splitting the augmented Lagrangian function $L(x_1,\ldots,x_p,\lambda)$ in the primal-dual order, and thus obtain the following scheme similar as (\ref{DADMM}):
\begin{equation}\label{PD-Prediction}
\left\{
\begin{array}{l}
 \tilde{x}_1^k  \in \arg\min \bigl\{ \theta_1(x_1)  -x_1^TA_1^T\lambda^k   +\frac{\beta}{2} \|A_1(x_1-x_1^k)\|^2    \;|\;   x_1\in{\cal X}_1  \bigr\};\\[0.2cm]
 \tilde{x}_2^k \in \arg\min \bigl\{\theta_2(x_2)  -x_2^TA_2^T\lambda^k   +\frac{\beta}{2} \|A_1(\tilde{x}_1^k-x_1^k)
                 + A_2(x_2-x_2^k)\|^2  \;|\;  x_2\in{\cal X}_2  \bigr\};\\[0.1cm]
  \qquad \qquad \vdots \\
  \tilde{x}_i^k\in\arg\min\bigl\{\theta_i(x_i)  -x_i^TA_i^T\lambda^k   +\frac{\beta}{2} \| \sum_{j=1}^{i-1}A_j(\tilde{x}_j^k-x_j^k)
                 + A_i(x_i-x_i^k)\|^2  \;|\;  x_i\in{\cal X}_i  \bigr\};\\[0.1cm]
  \qquad  \qquad \vdots   \\
  \tilde{x}_{p}^k \in \arg\min \bigl\{
    \theta_{p}(x_{p})  -x_{p}^TA_{p}^T\lambda^k   +\frac{\beta}{2} \| \sum_{j=1}^{p-1}A_j(\tilde{x}_j^k-x_j^k)
                 + A_{p}(x_{p}-x_{p}^k)\|^2
  \;|\;  x_p\in{\cal X}_{p}  \bigr\}; \\[0.3cm]
    \tilde{\lambda}^k =\arg\max\bigl\{- \lambda^T\bigl(\sum_{j=1}^{p} A_j\tilde{x}_j^k -b\bigr) -
       \frac{1}{2\beta}\|\lambda-\lambda^k\|^2   \;|\;  \lambda\in \Lambda \bigr\}.
     \end{array}
\right.
\end{equation}
Note that linear inequalities are considered in the model (\ref{Problem-m}). Hence, updating the Lagrange multiplier $\lambda$ in (\ref{PD-Prediction}) involves the constraint $\Lambda$. For the same reason in \cite{CHYY},  the coarse splitting scheme (\ref{PD-Prediction}) has no guaranteed convergence, but it can be used as a starting point to design a splitting contraction algorithm with provable convergence. As analyzed in \cite{HeXuYuan}, the coarse splitting scheme (\ref{PD-Prediction}) can be rewritten as
\begin{subequations}  \label{xPD-PreD}
  \[\label{xPD-PreD-F }  \tilde{w}^k\in \Omega,\quad
    \theta(x) - \theta(\tilde{x}^k) +(w- \tilde{w}^k)
   ^T F(\tilde{w}^k) \ge (w- \tilde{w}^k)^T
    Q_{\pd}(w^k - \tilde{w}^k),
     \quad \forall \, w\in {\Omega},
       \]
where
\[  \label{xPD-Q}
  Q_{\pd} = \left(\begin{array}{ccccc}
         \beta A_1^TA_1   &      0                         & \cdots      &  0           &    A_1^T   \\[0.3cm]
         \beta A_2^TA_1  &   \beta A_2^TA_2 &  \ddots    &  \vdots   &  A_2^T    \\[0.3cm]
             \vdots              &  &   \ddots          &       0        &    \vdots \\[0.3cm]
         \beta A_{p}^TA_1  &   \beta A_{p}^TA_2      &  \cdots  &\beta A_{p}^TA_{p}  &   A_{p}^T\\[0.3cm]
               0  &     0&  \cdots &   0   &  \frac{1}{\beta}I_{m}
 \end{array}\!\!\right), \]
 \end{subequations}
which is in form of the prediction step (\ref{Frame-Q}) with the prediction matrix $Q_{\pd}$. It is easy to see that the matrix $Q_{\pd}^T+Q_{\pd}$ is positive definite if all $A_i$'s are full column rank.

To simplify the notation, let us further denote
\[ \label{xPD-Bw}
   P =
       \left(\begin{array}{ccccc}
         \sqrt{\beta} A_1    &   \qquad 0 \qquad & \qquad \cdots \qquad & \qquad \cdots \qquad  &    0   \\[0.1cm]
         0   &   \sqrt{\beta} A_2   &  \ddots   &     &    \vdots \\[0.1cm]
             \vdots    &  \ddots   &   \ddots          &    \ddots             &    \vdots \\[0.1cm]
          \vdots   &            &  \ddots  &\sqrt{\beta} A_{p}   &   0\\[0.1cm]
              0   &    \cdots &  \cdots &    0   &  \frac{1}{\sqrt{\beta}}I_{m}
 \end{array}\!\!\right),   \qquad   \xi= Pw= \left(\begin{array}{c}
          \sqrt{\beta} A_1x_1  \\[0.1cm]
             \sqrt{\beta} A_2x_2   \\[0.1cm]
             \vdots \\[0.1cm]
               \sqrt{\beta} A_px_p \\[0.1cm]
            \frac{1}{\sqrt{\beta}} \lambda
   \end{array} \right);\]
and also
$${\Xi}= \big\{\xi  \;|\;  \xi= Pw, \; w\in \Omega\big\}\quad   \hbox{and}  \quad   {\Xi^*}= \big\{\xi^*  \;|\;  \xi^*= Pw^*, \; w^*\in \Omega^*\big\}. $$
Then, with $P$ and $\xi$ in \eqref{xPD-Bw}, we can rewrite the VI (\ref{xPD-PreD-F })-(\ref{xPD-Q}) as
\[   \label{xPD-wkw}
     \theta(x) - \theta(\tilde{x}^k)  +  (w - \tilde{w}^k)^T F(\tilde{w}^k) \ge
      (\xi  -\tilde{\xi}^k)^T{\Q}_{\pd}(\xi^k-\tilde{\xi}^k), \quad \forall  \, w \in {\Omega},
  \]
where
 \[  \label{xMatrixBQ-PD}  Q_{\pd}= P^T {\Q}_{\pd} P \quad \hbox{with}\quad {\Q}_{\pd} = \left(\begin{array}{ccccc}
             I_m   &      0      & \cdots      &  0           &     I_m \\[0.1cm]
             I_m   &     I_m  &  \ddots    &  \vdots   &     I_m   \\[0.1cm]
             \vdots   &  &   \ddots          &       0        &    \vdots \\[0.1cm]
             I_m   &     I_m      &  \cdots  &   I_m &      I_m\\[0.1cm]
               0  &     0&  \cdots &   0   &   I_m
 \end{array}\!\!\right).
  \]

\subsubsection{Dual-primal order}  \label{Subsec-mDP}

In (\ref{PD-Prediction}), the primal variables are solved first, followed by the dual variable. We can alternatively update the dual variable first. That is, we can consider the following scheme which updates the variables in the dual-primal order:

\begin{equation}\label{DP-Prediction}
\left\{
\begin{array}{l}
 \tilde{\lambda}^k =\arg\max\bigl\{- \lambda^T\bigl(\sum_{j=1}^{p} A_j{x}_j^k -b\bigr) - \frac{1}{2\beta}\|\lambda-\lambda^k\|^2    \;|\;  \lambda\in \Lambda \bigr\}; \\[0.2cm]
 \tilde{x}_1^k  \in \arg\min \bigl\{ \theta_1(x_1)  -x_1^TA_1^T\tilde{\lambda}^k   +\frac{\beta}{2} \|A_1(x_1-x_1^k)\|^2  \;|\;   x_1\in{\cal X}_1  \bigr\};  \\[0.2cm]
 \tilde{x}_2^k \in \arg\min \bigl\{\theta_2(x_2)  -x_2^TA_2^T\tilde{\lambda}^k   +\frac{\beta}{2} \|A_1(\tilde{x}_1^k-x_1^k) + A_2(x_2-x_2^k)\|^2  \;|\;  x_2\in{\cal X}_2  \bigr\};  \\[0.1cm]
  \qquad \qquad \vdots \\
  \tilde{x}_i^k\in \arg\min\bigl\{\theta_i(x_i)  -x_i^TA_i^T\tilde{\lambda}^k   +\frac{\beta}{2} \| \sum_{j=1}^{i-1}A_j(\tilde{x}_j^k-x_j^k)
                 + A_i(x_i-x_i^k)\|^2  \;|\;  x_i\in{\cal X}_i  \bigr\};\\[0.1cm]
  \qquad  \qquad \vdots   \\
  \tilde{x}_{p}^k \in\arg\min \bigl\{
    \theta_{p}(x_{p})  -x_{p}^TA_{p}^T\tilde{\lambda}^k   +\frac{\beta}{2} \| \sum_{j=1}^{p-1}A_j(\tilde{x}_j^k-x_j^k)
                 + A_{p}(x_{p}-x_{p}^k)\|^2   \;|\;  x_p\in{\cal X}_{p}  \bigr\} .
     \end{array}
\right.
\end{equation}

Similarly as (\ref{PD-Prediction}), the coarse splitting scheme (\ref{DP-Prediction}) has no guaranteed convergence, but it can be used as a starting point to design a splitting contraction algorithm with provable convergece. As analyzed in \cite{HeXuYuan}, the scheme (\ref{DP-Prediction}) can be rewritten as
\begin{subequations}  \label{yDP-PreD}
  \[\label{yDP-PreD-F }  \tilde{w}^k\in \Omega,\quad
    \theta(x) - \theta(\tilde{x}^k) +(w- \tilde{w}^k)
   ^T F(\tilde{w}^k) \ge (w- \tilde{w}^k)^T
    Q_{\dup}(w^k - \tilde{w}^k),
     \; \;\forall\; w\in {\Omega},
       \]
where
\[  \label{yDP-Q}
  Q_{\dup} = \left(\begin{array}{ccccc}
         \beta A_1^TA_1   &      0                         & \cdots      &  0           &   0   \\[0.3cm]
         \beta A_2^TA_1  &   \beta A_2^TA_2 &  \ddots    &  \vdots   &  0   \\[0.3cm]
             \vdots              &  &   \ddots          &       0        &    \vdots \\[0.3cm]
         \beta A_{p}^TA_1  &   \beta A_{p}^TA_2      &  \cdots  &\beta A_{p}^TA_{p}  &  0\\[0.3cm]
              - A_1  &   -  A_2&  \cdots &  - A_p   &  \frac{1}{\beta}I_{m}
 \end{array}\!\!\right). \]
 \end{subequations} It is easy to see that the matrix $Q_{\dup}^T+Q_{\dup}$ is positive definite if all $A_i$'s are full column rank.

Also, using $P$ and $\xi$  in \eqref{xPD-Bw},  we can rewrite the VI (\ref{yDP-PreD-F })-(\ref{yDP-Q}) as
\[   \label{yDP-wkw}  \tilde{w}^k\in \Omega,\quad
     \theta(x) - \theta(\tilde{x}^k)  +  (w - \tilde{w}^k)^T F(\tilde{w}^k) \ge
      (\xi  -\tilde{\xi}^k)^T{\Q}_{\dup}(\xi^k-\tilde{\xi}^k), \quad \forall  \, w \in {\Omega}.
  \]
where
   \[  \label{yMatrixBQ-DP}   Q_{\dup}= P^T {\Q}_{\dup} P\quad\hbox{with}\quad \Q_{\dup}=  \left(\begin{array}{ccccc}
             I_m   &      0      & \cdots      &  0           &     0 \\[0.1cm]
             I_m   &     I_m  &  \ddots    &  \vdots   &     0   \\[0.1cm]
             \vdots   &  &   \ddots          &       0        &    \vdots \\[0.1cm]
             I_m   &     I_m      &  \cdots  &   I_m &      0\\[0.1cm]
             - I_m  &  - I_m &  \cdots &  - I_m  &   I_m
 \end{array}\!\!\right).
  \]

%
%
%
%
%
%

\subsection{Representation of the prediction-correction framework (\ref{Framework}) and convergence conditions (\ref{Frame-HG})}

To show splitting contraction algorithms for the model (\ref{Problem-m}) more clearly, we can rewrite the prediction-correction framework (\ref{Framework}) and convergence conditions (\ref{Frame-HG}) in the context of the VI (\ref{VI-pFORM-Q})-(\ref{VI-pFORM-F}) with the notation $P$ and $\xi$ in (\ref{xPD-Bw}).

\begin{center}\fbox{\begin{minipage}{15.5cm}

\begin{subequations} \label{m-PRECOR}
\noindent{\bf \normalsize Prediction-correction framework for the VI (\ref{VI-pFORM-Q})-(\ref{VI-pFORM-F})}

[Prediction Step.]  With given $\xi^k= Pw^k$, find $\tilde{w}^k \in \Omega$ such that
 \[  \label{m-PRE} \tilde{w}^k \in \Omega, \;\;\theta(x) - \theta(\tilde{x}^k)  +  (w - \tilde{w}^k)^T F(\tilde{w}^k) \ge
      (\xi  -\tilde{\xi}^k)^T{\Q}(\xi^k-\tilde{\xi}^k), \;\; \forall  \, w \in {\Omega},
        \]
 where the matrix $\Q \in \Re^{(p+1)m\times (p+1)m}$ is not necessarily symmetric but
 the matrix ${\Q}^T+{\Q}$  is assumed to be positive definite.

\smallskip
[Correction Step.] Find a nonsingular matrix $\M$ and update $\xi$ by
 \[ \label{m-COR}   {\xi}^{k+1} = {\xi}^k -  \M(\xi^k - \tilde{\xi}^k).\]
\end{subequations}

\end{minipage}} \end{center}

\begin{center}\fbox{\begin{minipage}{15.5cm}
\noindent{\bf \normalsize Convergence conditions}

For the matrices $Q$ and $M$ used in (\ref{m-PRE}) and (\ref{m-COR}), respectively, there exists a matrix $\mathcal{H}\succ0$ such that
 \begin{subequations} \label{m-Frame-HG}
\[   \label{m-Frame-H} \mathcal{H}\M=\Q,\]
and
\[   \label{m-Frame-G}
\mathcal{G}:=  \Q^T +\Q - \M^T\mathcal{H}\M \succ 0.
       \]
\end{subequations}
\end{minipage}
}
\end{center}

Obviously, the VIs \eqref{xPD-wkw} and \eqref{yDP-wkw} are specific cases of (\ref{m-PRE}), whose corresponding prediction matrices are $\Q_{\pd}$ in (\ref{xMatrixBQ-PD}) and $\Q_{\dup}$ in (\ref{yMatrixBQ-DP}), respectively.

\subsection{Specifications of splitting contraction algorithms}  \label{Sec-Correction}

Now, we focus on specifying the correction matrix $\M$ in (\ref{m-COR}) with $\Q_{\pd}$ in (\ref{xMatrixBQ-PD}) and $\Q_{\dup}$ in (\ref{yMatrixBQ-DP}). We will show that the algorithms in \cite{HeXuYuan} can be recovered, while more new algorithms can be designed easily by our proposed construction strategies. With $\Q_{\pd}$ and $\Q_{\dup}$, their associated correction matrices are denoted by $\M_{\pd}$ and $\M_{\dup}$, and $\G_{\pd}$ and $\G_{\dup}$, respectively.


\subsubsection{Some matrices}

To further simplify the notation to be used, we define the following $p\times p$ block matrices:
   \[  \label{x-LowL}
         {\cal L} =  \left(\begin{array}{cccc}
             I_m   &      0      & \cdots      &  0           \\[0.1cm]
             I_m   &     I_m  &  \ddots    &  \vdots   \\[0.1cm]
             \vdots   &  &   \ddots          &       0       \\[0.1cm]
             I_m   &     I_m      &  \cdots  &   I_m
         \end{array}\!\!\right)\qquad \;\hbox{and}\;   \qquad    {\cal I} =
  \left(\begin{array}{cccc}
             I_m   &      0      & \cdots      &  0           \\[0.1cm]
                  0   &     I_m  &  \ddots    &  \vdots   \\[0.1cm]
             \vdots   & \ddots &   \ddots          &       0       \\[0.1cm]
               0  &     \cdots     &  0 &   I_m
         \end{array}\!\!\right).
             \]
We also define the $1\times p$  block matrix
\[  \label{x-Row-E}     {\cal E} =  \left(\!\begin{array}{cccc}
                        I_m   &     I_m      &  \cdots  &   I_m
         \end{array}\!\!\right).   \]
It is cleat that
\[  \label{LTL}   {\cal L}^T + {\cal L} = {\cal I}  + {\cal E}^T {\cal E}.     \]
Furthermore, the matrix ${\Q}_{\pd}$ in \eqref{xMatrixBQ-PD}  has the form
$$  {\Q}_{\pd} =\left(\begin{array}{cc}
                         {\cal  L} &      {\cal E}^T\\[0.1cm]
               0  &   I_m
 \end{array}\!\!\right)\qquad   \hbox{and thus} \qquad   {\Q}^T_{\pd} +   {\Q}_{\pd}   = \left(\begin{array}{cc}
                          {\cal I}  + {\cal E}^T {\cal E} &   {\cal E}^T \\[0.1cm]
                   {\cal E}  &  2 I_m
 \end{array}\!\!\right).   $$
 Similarly,  the matrix ${\Q}_{\dup}$ in \eqref{yMatrixBQ-DP}  has the form
 $$    {\Q}_{\dup}   = \left(\begin{array}{cc}
                         {\cal  L} &   0  \\[0.1cm]
                   - {\cal E}  &   I_m
 \end{array}\!\!\right)   \qquad   \hbox{and thus} \qquad   {\Q}^T_{\dup} +   {\Q}_{\dup}   = \left(\begin{array}{cc}
                          {\cal I}  + {\cal E}^T {\cal E} &   -{\cal E}^T \\[0.1cm]
                  - {\cal E}  &   2  I_m
 \end{array}\!\!\right).    $$

To further analyze the correction steps associated with the correction matrices $\M_{\pd}$ and $\M_{\dup}$, let us take a closer look at  the matrices ${\Q}_{\pd}^{-T}$  and ${\Q}_{\dup}^{-T}$. Indeed, we have
\[  \label{CQ-Tpd}    {\Q}_{\pd}^{-T} =  \left(\begin{array}{cc}
                         {\cal  L}^T &     0  \\[0.1cm]
                 {\cal E}   &   I_m
 \end{array}\!\!\right)^{-1} =     \left(\begin{array}{cc}
                         {\cal  L}^{-T} &  0   \\[0.1cm]
               -  {\cal E} {\cal  L}^{-T}   &   I_m
 \end{array}\!\!\right).
\]
and
\[   \label{CQ-Tdp}   {\Q}_{\dup}^{-T} =  \left(\begin{array}{cc}
                         {\cal  L}^T &    - {\cal E}^T  \\[0.1cm]
                 0    &   I_m
 \end{array}\!\!\right)^{-1} =     \left(\begin{array}{cc}
                         {\cal  L}^{-T} &     {\cal  L}^{-T}  {\cal E}^T  \\[0.1cm]
                 0   &   I_m
 \end{array}\!\!\right).
\]
Recall the respective definitions ${\cal L}$ and ${\cal E}$ in \eqref{x-LowL} and \eqref{x-Row-E}. We have
\[\label{LT-1}
           {\cal L}^{-T} =  \left(\begin{array}{cccc}
             I_m   &   - I_m     &    0    &         0 \\[0.1cm]
               0  &     I_m  &  \ddots    &            0  \\[0.1cm]
             \vdots   &  \ddots  &   \ddots          &       - I_m    \\[0.1cm]
              0   &    \cdots     &  0  &  I_m
               \end{array}\!\!\right) \]
and
\[\label{LT-2}
  {\cal E} {\cal L}^{-T}
                  = \left(\!\begin{array}{cccc}
                        I_m  \! &\!    0   \!   & \! \cdots \! &  0
         \end{array}\!\!\right)\quad\hbox{and}\quad   {\cal  L}^{-T}  {\cal E}^T =
           \left(\begin{array}{c}
                   0 \\[0.1cm]
                   \vdots  \\[0.1cm]
               0    \\[0.1cm]
              I_m
               \end{array}\!\!\right). \]
Hence, the matrices $ {\Q}_{\pd}^{-T}$ in (\ref{CQ-Tpd}) and ${\Q}_{\dup}^{-T}$ in (\ref{CQ-Tdp}) are both very simple in structure; their entries only consist of blocks of $I_m$, $-I_m$ and $0$.

\subsubsection{$\M_{\pd}$ for the primal-dual prediction (\ref{PD-Prediction})}

With (\ref{PD-Prediction}) as the prediction step, the prediction matrix ${\Q}_{\pd}$ is given in (\ref{xMatrixBQ-PD}). To construct the corresponding correction matrix $\M_{\pd}$, for example, we can choose
\[  \label{PD-D}  {\mathcal{D}_{\pd}}  =
 \left(\begin{array}{cc}
        \nu{\cal I}  &0   \\[0.1cm]
                 0 &  I_m
 \end{array}\!\!\right),\]
with any $\nu \in (0,1)$. Recall (\ref{DG-QTQ}). Thus, (\ref{PD-D}) also means
 \begin{eqnarray*}
   {\G}_{\pd}  &:= & (  {\Q}_{\pd}^T +  {\Q}_{\pd})   -{\mathcal{D}_{\pd}}   =  \left(\begin{array}{cc}
        (1-\nu){\cal I} + {\cal E}^T {\cal E}  &      {\cal E}^T\\[0.1cm]
                 {\cal E}  &  I_m
 \end{array}\!\!\right).
\end{eqnarray*}
It is clear that both the matrices ${\mathcal{D}}_{\pd}$  and  ${\G}_{\pd}$ are positive definite. According to (\ref{E-aDM}), (\ref{CQ-Tpd}) and (\ref{PD-D}), the correction matrix ${\M}_{\pd}$ can be constructed as
\[   \label{PD-M}   {\M}_{\pd} =  {\Q}_{\pd}^{-T}{\mathcal{D}_{\pd}} =\left(\begin{array}{cc}
                         {\cal  L}^{-T} &  0   \\[0.1cm]
               -  {\cal E} {\cal  L}^{-T}   &   I_m
 \end{array}\!\!\right) \left(\begin{array}{cc}
        \nu{\cal I}  &0   \\[0.1cm]
                 0 &  I_m
 \end{array}\!\!\right) =
    \left(\begin{array}{cc}
        {\nu}{\cal L}^{-T}  &  0\\[0.2cm]
     -\nu {\cal E} {\cal L}^{-T}  &   I_m
        \end{array}\right).
   \]
This coincides with the correction step in Section 7 of \cite{HeXuYuan}. Recall (\ref{LT-1}) and (\ref{LT-2}). We know that the correction step (\ref{m-COR}) with the correction matrix ${\M}_{\pd}$ defined in (\ref{PD-M}) is extremely easy to be implemented. Hence, the implementation of the specified splitting contraction algorithm mainly needs to solve the splitting $x_i$-subproblems in (\ref{PD-Prediction}).
%

\subsubsection{$\M_{\dup}$ for the dual-primal  prediction (\ref{DP-Prediction})}

With (\ref{DP-Prediction}) as the prediction step, the prediction matrix ${\Q}_{\dup}$ is given in (\ref{yMatrixBQ-DP}). To construct the  correction matrix $\M_{\dup}$, for example, we can choose
  \[   \label{DP-D}   {\cal{D}_{\dup}}
   =  \left(\begin{array}{cc}
            {\nu}{\cal I} + {\cal E}^T{\cal E}  &    -{\cal E}^T \\[0.1cm]
              - {\cal E}  &   I_{m}
        \end{array}\right)  \]
with $\nu \in(0,1)$. Recall (\ref{DG-QTQ}). Thus, (\ref{DP-D}) also means
 \begin{eqnarray*}
   {\G}_{\dup}  &:= & (  {\Q}_{\dup}^T +  {\Q}_{\dup})   -{\mathcal{D}_{\dup}}   =  \left(\begin{array}{cc}
        \nu{\cal E}  &  0\\[0.1cm]
                   0  &  I_m
 \end{array}\!\!\right).
 \end{eqnarray*}
It is clear that both the matrices ${\mathcal{D}}_{\dup}$  and  ${\G}_{\dup}$ are positive definite.  According to (\ref{E-aDM}), (\ref{CQ-Tdp}) and (\ref{DP-D}), the correction matrix ${\M}_{\dup}$ can be constructed by
\[   \label{DP-M}   {\M}_{\dup} =  {\Q}_{\dup}^{-T}{\mathcal{D}_{\dup}} =\left(\begin{array}{cc}
                         {\cal  L}^{-T} &    {\cal  L}^{-T}{\cal E}^T   \\[0.1cm]
                                        0    &   I_m
 \end{array}\!\!\right)
    \left(\begin{array}{cc}
            {\nu}{\cal I} + {\cal E}^T{\cal E}  &    -{\cal E}^T \\[0.1cm]
              - {\cal E}  &   I_{m}
        \end{array}\right)  =  \left(\begin{array}{cc}
        {\nu}{\cal L}^{-T}  &  0\\[0.2cm]
     - {\cal E} &   I_m
        \end{array}\right).
   \]
This coincides with the correction step in Section 8 of \cite{HeXuYuan}.
Recall (\ref{LT-1}) and (\ref{LT-2}). We know that the correction step (\ref{m-COR}) with the correction matrix ${\M}_{\dup}$ defined in (\ref{DP-M}) is also extremely easy to be implemented. Hence, the implementation of the specified splitting contraction algorithm mainly needs to solve the splitting $x_i$-subproblems in (\ref{DP-Prediction}).

\subsubsection{More choices}

To specify the prediction-correction framework (\ref{m-PRECOR}) and ensure the convergence conditions (\ref{m-Frame-HG}) with a given $\Q$ satisfying $\Q^T+\Q \succ 0$, like (\ref{DG-QTQ}), the matrices $\mathcal{D}$ and $\mathcal{G}$ can be chosen with the only restriction
 \[\label{m-DG-QTQ}   \mathcal{D}\succ 0, \quad \mathcal{G}\succ 0, \quad \hbox{and}\quad   \mathcal{D} + \mathcal{G} = \Q^T+\Q. \]
Hence, there are infinitely many ways to construct $\mathcal{H}$ and $\M$ with the given prediction matrix ${\Q}_{\pd}$ in (\ref{xMatrixBQ-PD}) or ${\Q}_{\dup}$ in (\ref{yMatrixBQ-DP}). For instance, we can choose $\mathcal{G}_{\pd}$ and $\mathcal{G}_{\dup}$, instead of $\mathcal{D}_{\pd}$ and $\mathcal{D}_{\dup}$, as the matrices defined in (\ref{PD-D}) and (\ref{DP-D}), respectively; or we can choose
$$
 {\mathcal{D}}= \alpha\bigl[\Q^T+\Q\bigr]  \;\; \hbox{and}\;\; \G= (1-\alpha)\bigl[\Q^T + \Q\bigr], \;\; \alpha\in(0,1).
  $$
All these choices lead to new splitting contraction algorithms for the model (\ref{Problem-m}) with provable convergence that are not covered in \cite{HeXuYuan}.

Finally, we would emphasize that, as shown in (\ref{CQ-Tpd})-(\ref{LT-2}), $\Q_{\pd}^{-T}$ and $\Q_{\dup}^{-T}$ are both very simple in structure, and they do not cause too much additional computation for constructing the correction matrix via ${\M}_{\pd} =  {\Q}_{\pd}^{-T}{\mathcal{D}_{\pd}}$ (see (\ref{PD-M})) or  ${\M}_{\dup} =  {\Q}_{\dup}^{-T}{\mathcal{D}_{\dup}} $ (see (\ref{DP-M})). This advantage makes it practical and adaptable to choose more application-tailored $\mathcal{D}_{\pd}$ and $\mathcal{D}_{\dup}$ for specific applications of the model (\ref{Problem-m}).


%

\section{Conclusions}
We revisited a unified framework for algorithmic design and convergence analysis that can capture a series of our previous works of desiging/analyzing splitting contraction algorithms for separable convex programming problems, and provided some construction strategies to specify this unified framework. By the proposed strategies, once a matrix (e.g., $D$ or $G$ as mentioned) is chosen, a splitting contraction algorithm with provable convergence can be automatically generated. There are many specification strategies, and the flexibility of choosing such a matrix enables us to design model-tailored/application-tailored splitting contraction algorithms with easy subproblems conveniently. We illustrated how to apply this construction principle to generate easily implementable ADMM-based algorithms for separable convex programming models with linear constraints. The same methodology can be applied to improve some other algorithms with theoretical or numerical disadvantages for other optimization problems.

\end{document}